\documentclass[10pt]{article}
\usepackage{multirow,bigdelim}  
\usepackage{booktabs}
\usepackage{topcapt}       
\usepackage{amssymb,amsmath,bm,graphicx}
\usepackage{epstopdf}   

\let\OLDthebibliography\thebibliography   
\renewcommand\thebibliography[1]{
  \OLDthebibliography{#1}
  \setlength{\parskip}{0pt}
  \setlength{\itemsep}{0pt plus 0.3ex}
}

\newtheorem{defn}{Definition}

\newtheorem{rmk}{Remark}
\newcommand{\R}{\mathbb{R}}
\newcommand{\bi}{\begin{itemize}}
\newcommand{\ei}{\end{itemize}}
\newcommand{\ben}{\begin{enumerate}}
\newcommand{\een}{\end{enumerate}}
\newcommand{\be}{\begin{equation}}
\newcommand{\ee}{\end{equation}}
\newcommand{\bea}{\begin{eqnarray}} 
\newcommand{\eea}{\end{eqnarray}}
\newcommand{\ba}{\begin{align}} 
\newcommand{\ea}{\end{align}}
\newcommand{\bse}{\begin{subequations}} 
\newcommand{\ese}{\end{subequations}}
\newcommand{\bc}{\begin{center}}
\newcommand{\ec}{\end{center}}
\newcommand{\bfi}{\begin{figure}}
\newcommand{\efi}{\end{figure}}

\newcommand{\tbox}[1]{{\mbox{\rm \tiny #1}}}
\newcommand{\mbf}[1]{{\mathbf #1}}

\newcommand{\vt}[2]{\left[\begin{array}{r}#1\\#2\end{array}\right]} 
\newcommand{\xx}{\mbf{x}}
\newcommand{\rr}{\mbf{r}}
\newcommand{\yy}{\mbf{y}}
\newcommand{\zz}{\mbf{z}}
\newcommand{\nn}{\mbf{n}}

\newcommand{\dd}{\mbf{d}}

\newcommand{\fref}[1]{Fig.~\ref{#1}}          
\newcommand{\sref}[1]{Sec.~\ref{#1}}          
\newcommand{\tref}[1]{Table~\ref{#1}}

\newcommand{\al}{\alpha}        
\newcommand{\eps}{\varepsilon}
\newcommand{\pro}{{P}}       
\newcommand{\NI}{{I}}        
\newcommand{\vep}{\varepsilon}
\newcommand{\ui}{u^\tbox{inc}}         
\newcommand{\ti}{{\theta^\tbox{inc}}}        
\newcommand{\Dr}{{\mathcal D}}    
\newcommand{\Sr}{{\mathcal S}}    
\newcommand{\sump}{{\sum_{p=1}^P}}        
\newcommand{\cc}[1]{c^{#1}}          
\newcommand{\cp}[1]{\cc{#1}_p}        
\newcommand{\cphi}[1]{{\cp{#1} \phi^{#1}_p}}        
\newcommand{\cphin}[2]{{\cp{#1} \frac{\partial \phi^{#1}_p}{\partial #2}}}  
\newcommand{\nrb}{K}       
\newcommand{\mpspack}{{\tt MPSpack}}
\newcommand{\bigO}{{\mathcal O}}
\newcommand{\Nden}{N_\tbox{den}}        
\newcommand{\Ntot}{{\cal N}}        

\begin{document}


\title{Robust fast direct integral equation solver for quasi-periodic scattering problems with a large number of layers}

\author{Min Hyung Cho\footnote{\small
Department of Mathematics, Dartmouth College, Hanover, NH 03755.
email: {\tt min.h.cho@dartmouth.edu}}
\;and
Alex H. Barnett\footnote{\small
Department of Mathematics, Dartmouth College, Hanover, NH 03755.
email: {\tt ahb@math.dartmouth.edu}}
}
\date{\today}
\maketitle


\begin{abstract}
\begin{sloppypar}  
We present a new boundary integral formulation
for time-harmonic wave diffraction from
two-dimensional structures with many layers
of arbitrary periodic shape, such as multilayer dielectric gratings
in TM polarization.
Our scheme is robust at all scattering parameters,
unlike the conventional quasi-periodic Green's function
method which fails whenever any of the layers
approaches a Wood anomaly.
%
%
We achieve this by a decomposition into
near- and far-field contributions.
The former uses the free-space Green's function in a
second-kind
integral equation on one period of the material interfaces and their
immediate left and right neighbors;
the latter uses
proxy point sources and small least-squares solves
(Schur complements) to represent the
remaining contribution from distant copies.
%
%
By using high-order discretization on interfaces (including those with corners),
the number of unknowns per layer is kept small.
We achieve overall linear complexity in the number of layers,
by direct solution of the resulting block tridiagonal system.
For device characterization we present an efficient method
to sweep over multiple incident angles, and show a $25\times$ speedup over
solving each angle independently.
We solve the scattering from a 1000-layer structure with $3\times 10^5$
unknowns
to 9-digit accuracy in 2.5 minutes on a 
desktop workstation.
\end{sloppypar}
\end{abstract}


\bfi  \centering  
(a)\raisebox{-3.6in}{\includegraphics[width=2.1in]{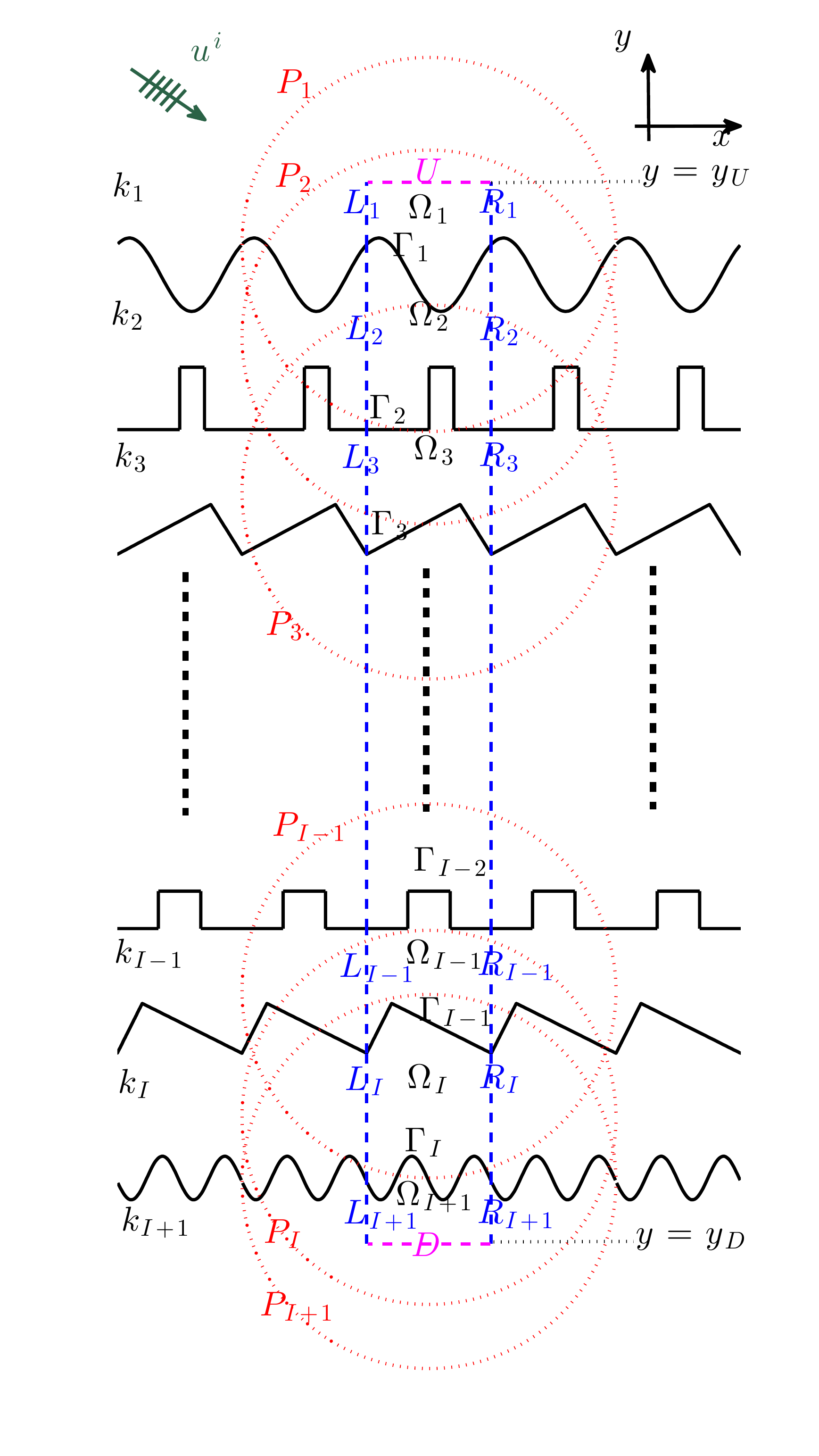}}
\hspace{.5in}
(b)\raisebox{-3.4in}{\includegraphics[width=2.1in]{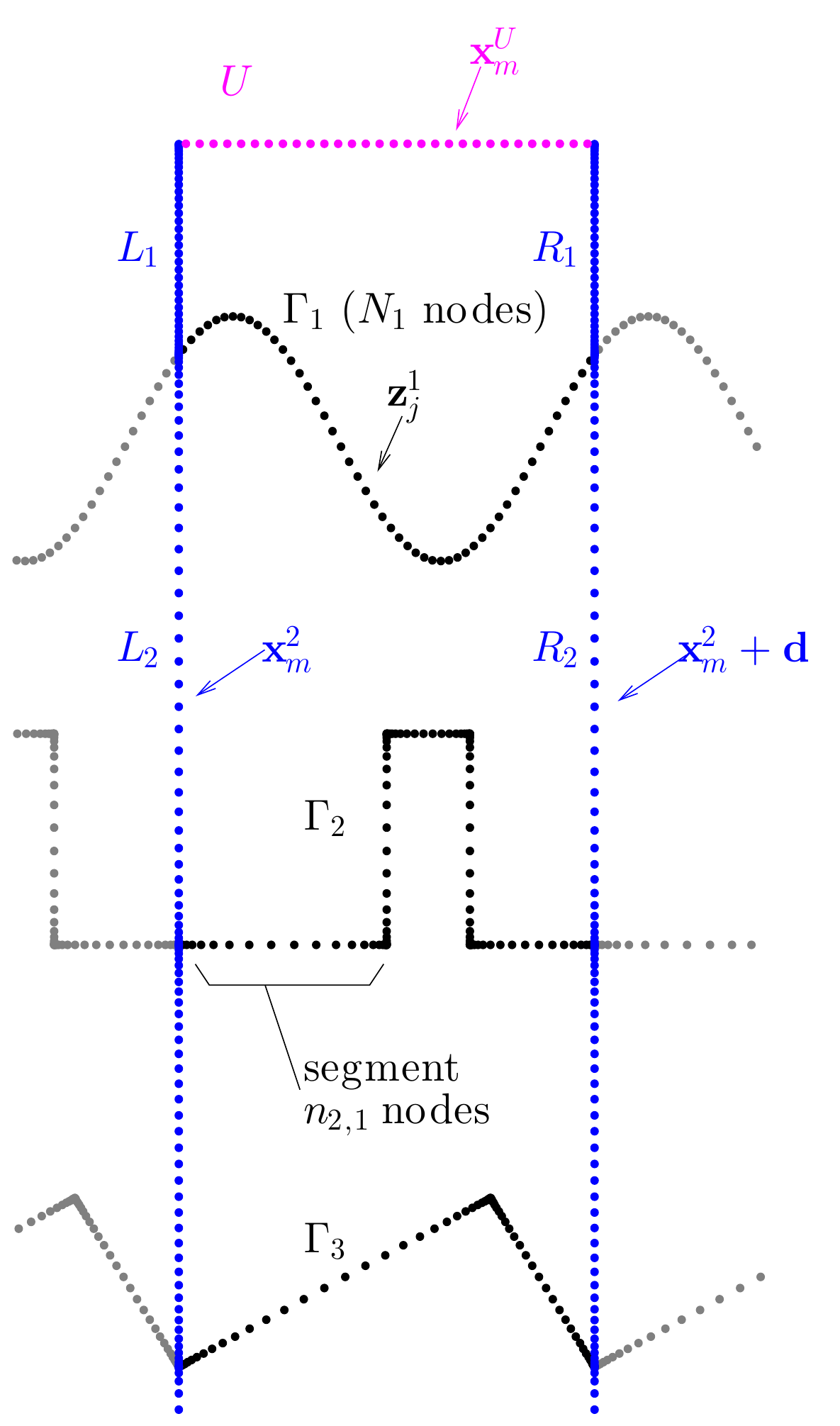}}
\caption{(a) Geometry of scattering problem. $\Gamma_i$ for $i=1,\ldots,I$ are the
    material interfaces with one spatial period lying between the vertical blue
dotted lines. The medium is uniform in the $i$th layer, which
lies above the $i$th interface and has wavenumber $k_i$.
Our algorithm also uses:
$L_i$ and $R_i$ which are the left and right walls of one period $\Omega_i$ of the
$i$th layer, $\pro_i$ the proxy circle for this layer,
and $U$, $D$ the upper and lower fictitious interfaces (at $y=y_U$ and $y=y_D$)
where the radiation condition is applied.
(b) Zoom of the top part of the geometry,
showing quadrature nodes for Nystr\"om method and
collocation (for clarity, less nodes are shown than actually used).
}\label{domain}
\efi

\section{Introduction}
Periodic geometries (such as diffraction gratings and antennae) and
multilayered media (such as dielectric mirrors)
are both essential for the manipulation of waves in
modern optical and electromagnetic devices.
In an increasing number of applications these features occur in tandem:
for instance, performance of a grating with transverse periodicity
will be enhanced by using many layers of different indices.
Dielectric gratings for high-powered laser
\cite{perryMLD,NIF} or wideband \cite{widebandgrating} applications rely on
structures with up to 50 dielectric layers.
Efficient thin-film photovoltaic cells exploit
multiple layers of silicon, transparent conductors, and dielectrics,
which are patterned to enhance absorption
\cite{atwater,kelzenberg}.
For solar thermal power,
efficient visible-light absorbers which reflect in the infrared
require patterning of several layers \cite{absorber}.
Related wave scattering problems appear in
photonic crystals \cite{jobook},
process control in semiconductor lithography \cite{scatterometry},
in the electromagnetic characterization of
increasingly multilayered integrated circuits,
and in models for underwater acoustic wave propagation.
%
In most of these settings, it is common to solve for the
scattering
for a large number of incident angles and/or wavelengths,
then repeat this inside a design optimization loop.
Thus, a robust and efficient solver is crucial.
We present such a solver, which scales optimally with respect to the
number of layers.

%
%

Let us describe the geometry of the problem (\fref{domain}(a)).
Consider $\NI$ interfaces $\Gamma_i$,
each of which has the same periodicity $d$ in the horizontal
($x$) direction.
The interfaces lie between $\NI+1$ homogeneous material layers,
each filling a domain $\Omega_i\subset \R^2$,
$i=1\,\ldots,\NI+1$. The $i$th layer lies between $\Gamma_i$ and $\Gamma_{i-1}$,
whilst the top and bottom layers are semi-infinite.
The wavenumber will be $k_i$ in the $i$th layer.
A plane wave is incident in the uppermost layer,
\be
\ui(\rr) = \left\{ \begin{array}{ll}e^{i\mathbf{k} \cdot \mathbf{r}},&  \rr\in\Omega_1\\
0, & \mbox{otherwise} \end{array} \right.
\label{ui}
\ee
with wavevector $\mathbf{k} = (k_1 \cos\ti, k_1 \sin\ti)$, at angle
$-\pi<\ti<0$.
The incident wave is {\em quasi-periodic} (periodic up to a phase),
meaning $\ui(x+d,y) = \al \ui(x,y)$ for all $(x,y) \in \R^2$,
where the {\em Bloch phase} (phase factor associated with translation by one
unit cell) is
\be
\al := e^{i d k_1 \cos \ti}
~.
\label{al}
\ee
Note that $\al$ is controlled by the period and indicent wave alone.
We will seek a solution sharing this quasi-periodic symmetry.

As is standard for scattering theory \cite{coltonkress},
the incident wave causes a scattered wave $u$ to be generated,
and the physical wave is their total $\ui + u$.
The scattered wave is given by solving the following
boundary value problem (BVP).
We have the Helmholtz equation in each layer,
\begin{equation}
\Delta u_i(\mathbf{r})+k_i^2 u_i(\mathbf{r}) = 0, \qquad \mathbf{r} = (x,y) \in \Omega_i
\label{layer_helmholtz}
\end{equation}
where we write $u_i$ for the scattered wave in the $i$th layer,
and the following interface, boundary, and radiation conditions:
\bi
\item Continuity of the value and derivative the total wave
on each interface, i.e.\
\begin{eqnarray}
u_1-u_2 &=& -\ui \qquad \mbox{ and } \quad \frac{\partial u_1}{\partial \mathbf{n}}-\frac{\partial u_2}{\partial \mathbf{n}} \;=\; -\frac{\partial \ui}{\partial \mathbf{n}} \qquad
\mbox{ on } \Gamma_1, \label{interface_cond1}\\
u_i-u_{i+1}& =& 0 \quad\qquad\; \mbox{ and } \quad \frac{\partial u_i}{\partial \mathbf{n}}-\frac{\partial u_{i+1}}{\partial \mathbf{n}} \;=\; 0 \qquad\quad \mbox{ on } \Gamma_i, \;i=2,3,\cdots I~.\label{interface_cond2}
\end{eqnarray}
\item Quasi-periodicity in all layers, i.e.\ for all $i=1,\ldots,I+1$,
\begin{equation}
u_i(x+d,y) = \alpha u_i(x,y)~, \qquad \mbox{ for all } (x,y) \in \R^2 ~.
\label{QP}
\end{equation}
%
%
\item Outgoing radiation conditions in $u_1$ and $u_{I+1}$,
namely the uniform convergence of Rayleigh--Bloch expansions \cite{bonnetBDS}
in the upper and lower half-spaces,
\begin{eqnarray}
u_1(x,y) &=& \sum_{n \in \mathbb{Z}} a^U_n e^{i\kappa_n x}e^{i k^U_n (y-y_U)}
\qquad \mbox{ for }  y\ge y_U \label{RB1}~,\\
u_{I+1}(x,y) &=& \sum_{n \in \mathbb{Z}} a^D_n e^{i\kappa_n x}e^{i k^D_n (-y+y_D)}
\qquad \mbox{ for }  y\le y_D\label{RB3}~,
\end{eqnarray}
where the horizontal wavenumbers in the modal expansion are
$$\kappa_n := k_1 \cos \ti+\frac{2\pi n}{d}~,$$
and
the upper and lower vertical wavenumbers are $k^U_n = \sqrt{k_1^2-\kappa_n^2}$
and $k^D_n = \sqrt{k_{I+1}^2-\kappa_n^2}$ respectively,
with all signs taken as positive real or imaginary.
The complex coefficients $a^U_n$ and $a^D_n$ are the Bragg diffraction amplitudes of the reflected and transmitted $n$th order; only the orders for which $k^U_n>0$ or $k^D_n>0$
carry power away from the system, and they characterize the
far-field diffraction pattern.
\ei

The above 
BVP describes time-harmonic
acoustics with layers of different sound speeds (but the same densities)
\cite{coltonkress}, or the time-harmonic
full Maxwell equations in the case of a $z$-invariant multilayer dielectric
in TM polarization, as we now briefly review;
see \cite{stratton, jackson}.
(The modification to the interface conditions for differing densities or TE polarization
are simple, and we will not present them.)
Letting the time dependence be $e^{-i\omega t}$, Maxwell's equations state that
the divergence-free vector fields $\mbf{E}$ and $\mbf{H}$ satisfy
\bea
\nabla \times \mathbf{E} &=& i\omega\mu \mathbf{H}~,\\
\nabla \times \mathbf{H} &=& -i\omega\varepsilon \mathbf{E}~.
\eea
Writing $\mathbf{E} = (0,0,u)$ and $\mathbf{H} = (H_x,H_y,0)$, and eliminating $\mbf{H}$
gives in the $i$th layer the Helmholtz PDE \eqref{layer_helmholtz} with wavenumber $k_i = \omega \sqrt{\vep_i \mu_i}$, where $\vep_i$ and $\mu_i$ are the permittivity and permeability 
of the layer.
Continuity of transverse $\mbf{E}$ and $\mbf{H}$ then gives the interface conditions
\eqref{interface_cond1}-\eqref{interface_cond2}.
Note that all of the layer wavenumbers $k_i$ are scaled linearly
by the overall frequency $\omega$.

The BVP \eqref{layer_helmholtz}--\eqref{RB3}
has a solution for all parameters \cite[Thm.~9.2]{bonnetBDS}.
The solution is unique at all but a discrete set of frequencies $\omega$ when
$\ti$ is fixed \cite[Thm.~9.4]{bonnetBDS};
these frequencies correspond to guided modes of the dielectric
structure, where resonance makes the physical problem ill-posed.
They are distinct from (but in the literature sometimes confused with)
{\em Wood anomalies} \cite{wood}, which are scattering parameters ($\ti$,$\omega$)
for which $k_U^n=0$ or $k_D^n=0$ for some $n$,
making the upper or lower
$n$th Rayleigh--Bloch mode a horizontally traveling plane wave.
A Wood anomaly does {\em not} prevent the solution from being unique,
although it does become arbitrarily sensitive to changes in $\ti$ or $\omega$.
For more detail see \cite{linton07,qpsc}, and the extensive review in the
three-dimensional (3D) case
by Shipman \cite{shipmanreview}.

There are many low-order numerical methods used to solve multilayer scattering
problems, which in test problems may only agree to 1 digit of accuracy \cite{refmodels}.
Finite difference time-domain (FDTD) \cite{taflove,qpfdtd}
is easy to code but has dispersion errors, and requires artificial absorbing
boundary conditions and arbitrarily long
settling times near resonances.
Direct discretization in the frequency-domain,
as in finite difference (FD) or finite element methods, is also possible
\cite{baodobsonrev,elschner,scatterometry},
although they require a large number of unknowns, and,
as the frequency grows, ``pollution'' error
means that an increasing number of degrees of freedom per wavelength is
needed \cite{pollution}.
They also demand artificial absorbing boundary conditions (perfectly matched layers).
The rigorous-coupled wave analysis (RCWA) or Fourier Modal Method
is specially designed for multilayer gratings
\cite{RCWAmoharam};
to overcome slow convergence a Fourier factorization method is needed
\cite{RCWAli,RCWAli4}.
However, RCWA is not easy to apply for arbitrary shapes (relying
on an intrinsically low-order ``staircase'' approximation of layer shapes),
nor to generalize to 3D.
Other methods include volume integral equations \cite{widebandgrating,lechleiter},
for which it is hard to exceed low-order convergence.
In general, when the layers are strictly planar, the problem becomes
1D \cite{chocai12} and the (unstable) transfer matrix and (stable) scattering matrix
approaches are standard \cite{transfer_matrix}.
However, we are concerned with interfaces of arbitrary shape.

Since the medium is piecewise constant, 
boundary integral equations (BIE) formulated on the interfaces are natural
and mathematically rigorous \cite{coltonkress,nedelec-starling,CWnystrom}.
By exploiting the reduced dimensionality,
the number of unknowns is much reduced,
and high-order quadratures exist, in 2D \cite{kress91,hao} but also in 3D.
Combined with fast algorithms for handling the resulting linear systems,
this leads to much higher efficiency and accuracy than FD methods, and
has started to be used effectively in periodic problems
\cite{CWnystrom,arenshabil,nicholas,otani08,brunohaslam09,qpfds,triplejuncQPFDS}.
In the setting of quasi-periodic scattering, the interfaces are unbounded
and the usual approach is to replace the 2D free-space kernel
(Green's function) for waves in layer $i$,
\be
G^i(\mbf{r},\mbf{r}') := \frac{i}{4} H^{(1)}_0(k_i\|\rr-\rr'\|)
\label{Gi}
\ee
(where $H_0^{(1)}$ is the Hankel function of order zero \cite{a+s}),
by the
quasi-periodic one obeying \eqref{QP}, in which case the problem may be formulated
on a {\em single period} of the interface.
The usual quasi-periodic Green's function for layer $i$ is
\be
G^i_{QP}(\mbf{r},\mbf{r}') := \sum_{l\in\mathbb{Z}} \al^l G^i(\mbf{r},\mbf{r}'+l\mbf{d})
~, \qquad \mbox{ where } \mbf{d}:=(d,0)~,
\label{GQPi}
\ee
a sum whose slow convergence renders it computationally useless.
Thus a large industry has been build around efficient evaluation of $G_{QP}$ using
convergence acceleration, Ewald's method, or lattice sums \cite{lintonrev}.
It is convenient to expand slightly the definition of Wood anomaly, as follows.
\begin{defn}  
We say that layer $i$ is at a Wood anomaly if either $\kappa_n = k_i$ or $\kappa_n = -k_i$ (or both) for some $n\in\mathbb{Z}$.
\end{defn}      
The problem then with $G_{QP}$-based methods
is that \eqref{GQPi} does not exist (the sum diverges)
whenever the $i$th layer is at a Wood anomaly.
As the number of different materials increases in a structure, the chances of
some layer hitting (or being close to)
a Wood anomaly, and thus of failure, increases. 
%

Two classes of solutions to this non-robustness problem have recently been introduced:
\ben
\item[I.]
Replace \eqref{GQPi} by the quasi-periodic Green's function for the
Dirichlet \cite{CWnystrom} or impedance half-plane problems \cite{horoshenkov},
or their generalization to larger numbers of images \cite{brunoqp3d}.
\item[II.]
Return to the free-space Green's function \eqref{Gi} for the central
unit cell plus immediate neighbors,
plus a new representation of far-field contributions,
imposing the quasi-periodicity condition in the least-squares sense via
additional rows in the linear system \cite{qplp,qpsc,qpfds}.
\een
In the multilayer dielectric setting robustness using class I would require
the impedance Green's function, which is difficult to evaluate;
while the class II contour-integral approach of Barnett--Greengard \cite{qpsc} does
not generalize well to multiple layers.

In this paper we introduce a simpler class II BIE method which combines the free-space
Green's function for the unit cell and neighbors, with a ring of 
{\em proxy point} sources (i.e.\ the method of fundamental solutions, or MFS \cite{Bo85,mfs})
to represent the far-field contributions which ``periodize'' the field.
This combines ideas in \cite[Sec.~3.2]{qplp} and the fast direct solver community
\cite[Sec.~5]{mdirect},
and has been independently proposed recently for Laplace problems
by Gumerov--Duraiswami \cite{gumerov}. Related representations have been used for
some time in the engineering community \cite{hafner90}.
A modern interpretation of the key idea is that the far-field contribution is
smooth---the interaction between distant periodic copies and the central unit cell is
{\em low rank}; eg see \cite{qpfds}---%
and hence only a small number of proxy points is needed, at least if the frequency is
not too high.
\begin{rmk} 
Conveniently, in our new formulation we can take Schur complements to
eliminate the proxy strength unknowns for each layer {\em without}
recreating \eqref{GQPi} and its associated Wood anomaly problem,
as would happen in prior methods \cite{qplp,qpsc} (see \cite[Remark~8]{qplp}).
The difference is that in \cite{qpsc} both upward and downward radiation conditions
(of the type \eqref{RB1} and \eqref{RB3})
are imposed on the Green's function, making it equivalent to \eqref{GQPi},
whereas we do not impose {\em any} radiation condition in the finite-thickness
layers, and impose outgoing conditions only in the semi-infinite layers 1 and $I+1$.
Since non-divergent Green's functions do exist which satisfy these minimal radiation conditions, they are selected by the least-squares linear algebra in the Schur complement
(see \sref{schur}).
\end{rmk}  
%
%
We present our new representation and its discretization in \sref{bie}, then combine
it in \sref{solve} with a direct solver which has two steps:
Schur complements to eliminate the proxy unknowns, followed by
direct block-tridiagonal factorization.
The tridiagonal structure arises simply because
layer $i$ couples only with
layers $i-1$ and $i+1$.
The overall scaling is $\bigO(IN^3)$, i.e.\ linear in the number of layers
and cubic in $N$ the number of unknowns per layer.
This allows our solver to tackle problems with $NI$ of order $10^6$ unknowns in
only a few minutes.
Since the solution is direct, as explained in \sref{sweep}
we can solve new incident waves that share the
same $\al$ without extra effort, and, by reusing matrix blocks, handle other
$\al$ values efficiently.
We test the solver's error and speed performance
with a variety of interface shapes, with up to $I=1000$ layers,
with random or periodic permittivities,
and multiple incident angles including a Wood anomaly,
in \sref{num}.
%
We conclude with a summary and implications for future work in \sref{conc}.

\section{Boundary integral formulation, periodizing scheme, and its
discretization} \label{bie}

We now reformulate the BVP as a system of linear second-kind
integral equations on the interfaces
$\Gamma_i$, $i=1,\ldots,I$ which lie in a single unit cell,
coupled with linear conditions on fictitious unit cell walls;
the complete system will be summarized by \eqref{system_bie} below.
A little extra geometry notation is needed, as shown in \fref{domain}.
Let us define the (central) unit cell as the vertical strip
of width $d$ lying between $x=-d/2$ and $x=d/2$; of course its
horizontal displacement is arbitrary.
The blue dashed vertical lines $\{L_i\}_{i=1}^{I+1}$ and $\{R_i\}_{i=1}^{I+1}$ are the left and right boundaries of the layer domains $\{\Omega_i\}_{i=1}^{I+1}$ lying inside the unit cell.
The proxy points for layer $i$ lie on the circle $\pro_i$ (shown by red dotted lines).
The magenta dashed lines $U$ and $D$ are fictitious interfaces for radiation conditions located at $y=y_U$ and $y=y_D$, touching $\Omega_1$ and $\Omega_{I+1}$, respectively.

\subsection{Representation of the scattered wave}
Using \eqref{Gi} we define standard potentials for the Helmholtz
equation, the single- and double-layer representations \cite{coltonkress}
lying on a general curve $W$, at wavenumber $k_i$ for the $i$th layer,
\be
(\Sr_{W}^i \sigma) (\mathbf{r}) := \int_{W} G^i(\mathbf{r}, \mathbf{r}') \sigma (\mathbf{r}') ~ds_{\mathbf{r}'}
~,
\quad
(\Dr_{W}^i \tau) (\mathbf{r}) := \int_{W} \frac{\partial G^i}{\partial \mathbf{n}'} (\mathbf{r}, \mathbf{r}') \tau (\mathbf{r}') ~ds_{\mathbf{r}'}
~,~\mbf{r}\in\R^2
\label{SD}
\ee
where $\mbf{n}'$ is the unit normal on the curve $W$ at $\mbf{r}'$,
and $ds$ the arclength element.
Shortly we will set $W$ to be either $\Gamma_{i-1}$ or $\Gamma_i$, with the normals
pointing down (into the layer below the interface).
Integral representations which include phased contributions from the nearest neighbors
are indicated with a tilde,
\be
(\tilde{\Sr}_{W}^i \sigma) (\mathbf{r}) :=
\sum_{l=-1}^1   \alpha^l \int_{W} G^i(\mathbf{r}, \mathbf{r}'+l\mathbf{d}) \sigma(\mathbf{r}') ~ds_{\mathbf{r}'}
~,
\qquad
(\tilde{\Dr}_{W}^i \tau) (\mathbf{r}) :=
\sum_{l=-1}^1  \alpha^l \int_{W}\frac{\partial G^i}{\partial \mathbf{n}'} (\mathbf{r}, \mathbf{r}'+l \mathbf{d}) \tau(\mathbf{r}') ~ds_{\mathbf{r}'}
~.
\label{SDsum}
\ee

Let the proxy points $\{\yy^i_p\}_{p=1}^P\in\R^2$ lie uniformly on the circle $\pro_i$
of radius $R$ which is centered on the domain $\Omega_i$.
As is well known
in MFS theory, increasing $R$ allows a higher convergence rate with respect to $P$
\cite[Thm.~3]{mfs}; however, since the proxy points are representing the contributions
from periodic interface copies $\{-\infty,\ldots,-3,-2\}$ and $\{2,3,\ldots,\infty\}$,
which thus have singularities at $|x|>3d/2$,
the proxy charge strengths will turn out to be
exponentially large \cite[Thm.~7]{mfs} if $R$ exceeds $3d/2$ by much.
In practice we choose $R\in[3d/2,2d]$.
In order to make the proxy representation more robust
at high wavenumbers we use a ``combined field'' approach,
choosing the proxy basis functions for $i$th layer,
\be
\phi^i_p(\rr) \;:=\; \frac{\partial G^i}{\partial \mathbf{n}_p} (\rr, \yy^i_p)
+ ik_i G^i(\rr,\yy^i_p)
~,\qquad \rr\in\Omega_i
~, \qquad p = 1,\ldots, P~
\label{phi}
\ee
where $\nn_p$ is the outwards-directed unit normal to the circle $P_i$
at the $p$th proxy point.
This results in smaller coefficients than if monopoles or dipoles alone were
used
(which can be justified by treating the proxy points
as a discrete approximation to a layer potential on $P_i$, and
considering arguments in \cite[Sec.~7.1]{zhaodet}).



Combining the near-field single- and double-layer potentials and proxy representations in each layer we have,
recalling the notation $u_i$ for $u$ in $\Omega_i$,
\begin{eqnarray}
u_1 &=&
\tilde{\Dr}^1_{\Gamma_1} \tau_1+ \tilde{\Sr}^1_{\Gamma_1} \sigma_1 +
\sump \cphi{1}
\label{firstlayer_sol}\\
u_i &=&
\tilde{\Dr}^i_{\Gamma_{i-1}} \tau_{i-1} + \tilde{\Sr}^i_{\Gamma_{i-1}} \sigma_{i-1} +
\tilde{\Dr}^i_{\Gamma_{i}} \tau_{i} + \tilde{\Sr}^i_{\Gamma_{i}} \sigma_{i} +
\sump \cphi{i}
~, \qquad i=2,3,\cdots, I
\label{middlelayer_sol}\\
u_{I+1} &=&
\tilde{\Dr}^{I+1}_{\Gamma_{I}} \tau_{I} + \tilde{\Sr}^{I+1}_{\Gamma_{I}} \sigma_{I} +
\sump \cphi{I+1}
\label{lastlayer_sol}
\end{eqnarray}
By construction,
for all layers $i=1,\ldots,I+1$,
and for all density functions $\sigma_i$ and $\tau_i$ and
proxy unknown vectors $\cc{i} := \{\cp{i}\}_{p=1}^P$,
this representation satisfies the Helmholtz equations \eqref{layer_helmholtz}.
In the following subsections, we describe in turn how to enforce the interface matching,
quasi-periodicity, and radiation conditions.
Each of these three conditions will comprise a block row of the final linear system
\eqref{system_bie} that enables us to solve for the densities and proxy unknowns.

\subsection{Matching conditions at material interfaces} \label{ic}

In this subsection, matching conditions \eqref{interface_cond1} and \eqref{interface_cond2} will be enforced at all material interfaces
in a standard M\"uller--Rokhlin \cite{muller,rokh83} scheme.

In the indirect approach, boundary integral operators arise from the restriction
of representations \eqref{SD} to curves \cite{coltonkress}.
Following \cite{qpsc} we use notation
$S^i_{V,W}$ to indicate the single-layer operator at wavenumber $k_i$
from a source curve $W$ to target curve $V$.
Similarly we use $D^i_{V,W}$ for the double-layer operator,
$D^{i,\ast}_{V,W}$ for the target-normal derivative of the single-layer operator,
and $T^i_{V,W}$ for the target-normal derivative of the double-layer operator.
As before, we use a tilde to indicate summation over the source curve and its
phased nearest neighbors, thus, for a target point $\xx\in V$,
\bea
(\tilde{S}_{V,W}^i \sigma) (\mathbf{r}) :=
\sum_{l=-1}^1  \! \alpha^l \!\int_{W} G^i(\mathbf{r}, \mathbf{r}'+l\mathbf{d}) \sigma(\mathbf{r}') ds_{\mathbf{r}'}
,
&&
\hspace{-2ex}
(\tilde{D}_{V,W}^i \tau) (\mathbf{r}) :=
\sum_{l=-1}^1 \! \alpha^l \!\int_{W}\frac{\partial G^i}{\partial \mathbf{n}'} (\mathbf{r}, \mathbf{r}'+l \mathbf{d}) \tau(\mathbf{r}') ds_{\mathbf{r}'}
\qquad\mbox{} 
\label{SDopsum}
\\
(\tilde{D}_{V,W}^{i,\ast} \sigma) (\mathbf{r}) :=
\sum_{l=-1}^1  \! \alpha^l \!\int_{W} \frac{\partial G^i}{\partial \nn}(\mathbf{r}, \mathbf{r}'+l\mathbf{d}) \sigma(\mathbf{r}') ds_{\mathbf{r}'}
,
&&
\hspace{-2ex}
(\tilde{T}_{V,W}^i \tau) (\mathbf{r}) :=
\sum_{l=-1}^1 \! \alpha^l \!\int_{W}\frac{\partial^2 G^i}{\partial \nn \partial \nn'} (\mathbf{r}, \mathbf{r}'+l \mathbf{d}) \tau(\mathbf{r}') ds_{\mathbf{r}'}
\qquad\mbox{} 
\label{DtTopsum}
\eea
When the target curve is the same as the source ($V=W$),
we note that the single-layer operator is a weakly singular integral operator,
that the action of the double-layer and its transpose must be taken in their
principal value sense, and that the $T$ operator is hypersingular.
%

At the first interface $\Gamma_1$, $u_1$ and $u_2$ are coupled.
The functions $u_1$, $u_2$, $\frac{\partial u_1}{\partial \mathbf{n}}$, and $\frac{\partial u_2}{\partial \mathbf{n}}$ at $\Gamma_1$ can be found by letting $\mathbf{r}$
in \eqref{firstlayer_sol}--\eqref{middlelayer_sol}
approach $\Gamma_1$
from the respective side, and using
the standard {\em jump relations} \cite[Thm.~3.1 and p.66]{coltonkress}
which introduce terms of one half times the density to each $D$ and $D^{\ast}$
term, giving
\begin{eqnarray}
u_1 &=&
-\frac{1}{2}\tau_1 + \tilde{D}^1_{\Gamma_1,\Gamma_1} \tau_1 + \tilde{S}^1_{\Gamma_1,\Gamma_1} \sigma_1 +
\sump \cphi{1}
~,\quad \mbox{ on } \Gamma_1
\\
u_2 &=&
\frac{1}{2}\tau_1 + \tilde{D}^2_{\Gamma_1,\Gamma_1} \tau_1+ \tilde{S}^2_{\Gamma_1,\Gamma_1} \sigma_1+\tilde{D}^2_{\Gamma_1,\Gamma_2} \tau_2+ \tilde{S}^2_{\Gamma_1,\Gamma_2} \sigma_2
+ \sump \cphi{2}
~,\quad \mbox{ on } \Gamma_1
\\
\frac{\partial u_1}{\partial \mathbf{n}} &=&
\tilde{T}^1_{\Gamma_1,\Gamma_1} \tau_1+\frac{1}{2}\sigma_1+ \tilde{D}^{1,*}_{\Gamma_1,\Gamma_1} \sigma_1 + \sump \cphin{1}{\nn}
~,\quad \mbox{ on } \Gamma_1
\\
\frac{\partial u_2}{\partial \mathbf{n}} &=&
\tilde{T}^2_{\Gamma_1,\Gamma_1} \tau_1-\frac{1}{2}\sigma_1+\tilde{D}^{2,*}_{\Gamma_1,\Gamma_1} \sigma_1+\tilde{T}^2_{\Gamma_1,\Gamma_2} \tau_2+ \tilde{D}^{2,*}_{\Gamma_1,\Gamma_2} \sigma_2
+ \sump \cphin{2}{\nn}
~,\quad \mbox{ on } \Gamma_1
~.
\end{eqnarray}

On this interface only, the matching conditions \eqref{interface_cond1} include
the indicent wave, and enforcing them using the above gives
the inhomogeneous coupled BIE and functional equations,
\begin{eqnarray}
-\tau_1+(\tilde{D}^1_{\Gamma_1,\Gamma_1}-\tilde{D}^2_{\Gamma_1,\Gamma_1})\tau_1+(\tilde{S}^1_{\Gamma_1,\Gamma_1}-\tilde{S}^2_{\Gamma_1,\Gamma_1})\sigma_1
-\tilde{D}^2_{\Gamma_1,\Gamma_2} \tau_2- \tilde{S}^2_{\Gamma_1,\Gamma_2} \sigma_2\nonumber\\
+ \;\sump (\cphi{1}-\cphi{2})|_{\Gamma_1} &=& -\ui|_{\Gamma_1} ~,
\qquad\mbox{}    
\label{bie_interface1}
\\
(\tilde{T}^1_{\Gamma_1,\Gamma_1}-\tilde{T}^2_{\Gamma_1,\Gamma_1})\tau_1+\sigma_1+(\tilde{D}^{1,*}_{\Gamma_1,\Gamma_1}-\tilde{D}^{2,*}_{\Gamma_1,\Gamma_1})\sigma_1
-\tilde{T}^2_{\Gamma_1,\Gamma_2} \tau_2 - \tilde{D}^{2,*}_{\Gamma_1,\Gamma_2} \sigma_2
\nonumber\\
+ \;
\sump \biggl.\biggl(\cphin{1}{\nn}-\cphin{2}{\nn}\biggr)\biggr|_{\Gamma_1} &=&
-\biggl.\frac{\partial \ui}{\partial \mathbf{n}}\biggr|_{\Gamma_1}
~.
\label{bie_interface2}
\eea
Note that the half density terms added to give a $-\tau_1$ in the first
equation and a $+\sigma_1$ in the second; these terms appear for every layer and
will cause the BIE to be of Fredholm second kind.

On the middle interfaces $\Gamma_i$, $i=2,\ldots,I-1$, we similarly match
$u_i$ and $u_{i+1}$ and their normal derivatives, noting that now there is
coupling to both the above and below interfaces,
but no effect of the incident wave, to get
\begin{eqnarray}
&&-\tau_i+(\tilde{D}^i_{\Gamma_i, \Gamma_i} - \tilde{D}^{i+1}_{\Gamma_i,\Gamma_i} )\tau_{i}
+ (\tilde{S}^i_{\Gamma_i, \Gamma_i}-\tilde{S}^{i+1}_{\Gamma_i,\Gamma_i}) \sigma_{i}+\tilde{D}^i_{\Gamma_i, \Gamma_{i-1}} \tau_{i-1}+ \tilde{S}^i_{\Gamma_i, \Gamma_{i-1}} \sigma_{i-1} \nonumber\\
&&-\tilde{D}^{i+1}_{\Gamma_i,\Gamma_{i+1}} \tau_{i+1}- \tilde{S}^{i+1}_{\Gamma_i,\Gamma_{i+1}} \sigma_{i+1}
+ \sump (\cphi{i}-\cphi{i+1})|_{\Gamma_i} = 0~,
\qquad\mbox{}    
\label{bie_interface3}
\\
&&(\tilde{T}^i_{\Gamma_i, \Gamma_i} - \tilde{T}^{i+1}_{\Gamma_i,\Gamma_i} )\tau_{i}+ \sigma_i
+ (\tilde{D}^{i,*}_{\Gamma_i, \Gamma_i}-\tilde{D}^{i+1,*}_{\Gamma_i,\Gamma_i}) \sigma_{i}+\tilde{T}^i_{\Gamma_i, \Gamma_{i-1}} \tau_{i-1}+ \tilde{D}^{i,*}_{\Gamma_i, \Gamma_{i-1}} \sigma_{i-1}\nonumber\\
&&-\tilde{T}^{i+1}_{\Gamma_i,\Gamma_{i+1}} \tau_{i+1}- \tilde{D}^{i+1,*}_{\Gamma_i,\Gamma_{i+1}} \sigma_{i+1}
+ \sump \biggl.\biggl(\cphin{i}{\nn}-\cphin{i+1}{\nn}\biggr)\biggr|_{\Gamma_i} = 0
~.
\label{bie_interface4}
\end{eqnarray}
On the bottom interface $\Gamma_I$, the only change
is the absence of coupling from any lower interface, so,
\begin{eqnarray}
-\tau_I+(\tilde{D}^I_{\Gamma_I,\Gamma_I}-\tilde{D}^{I+1}_{\Gamma_I,\Gamma_I})\tau_I+(\tilde{S}^I_{\Gamma_I,\Gamma_I}-\tilde{S}^{I+1}_{\Gamma_I,\Gamma_I})\sigma_I
+\tilde{D}^{I}_{\Gamma_I,\Gamma_{I-1}} \tau_{I-1}+ \tilde{S}^{I}_{\Gamma_I,\Gamma_{I-1}} \sigma_{I-1}
\nonumber\\
+ \; \sump (\cphi{I}-\cphi{I+1})|_{\Gamma_I} &=& 0~,
\label{bie_interface5}
\\
(\tilde{T}^I_{\Gamma_I,\Gamma_I}-\tilde{T}^{I+1}_{\Gamma_I,\Gamma_I})\tau_I+\sigma_I+(\tilde{D}^{I,*}_{\Gamma_I,\Gamma_I}-\tilde{D}^{I+1,*}_{\Gamma_I,\Gamma_I})\sigma_I
+\tilde{T}^I_{\Gamma_I,\Gamma_{I-1}} \tau_{I-1} + \tilde{D}^{I,*}_{\Gamma_{I},\Gamma_{I-1}} \sigma_{I-1}
\nonumber \\
+ \; \sump \biggl.\biggl(\cphin{I}{\nn}-\cphin{I+1}{\nn}\biggr)\biggr|_{\Gamma_I} &=& 0~.
\label{bie_interface6}
\end{eqnarray}

We wish to write these in a more compact form, hence we pair up double- and single-layer
densities, then stack them into a single column vector,
\begin{equation}
\bm{\eta} := \left[\begin{array}{c}\bm{\eta}_1, \bm{\eta}_2  , \cdots , \bm{\eta}_I \end{array}\right]^T
, \qquad \mbox{ where } \quad
\bm{\eta}_i := \vt{\tau_i}{\sigma_i}, \quad i=1,2,\cdots,I
~.
\ee
Similarly we stack the proxy strength vectors
$\cc{i} = \{\cp{i}\}_{p=1}^P$,
and form a vector of right-hand side functions,
\begin{equation}
\mathbf{c} = \left[\begin{array}{c}c^1 , c^2 , \cdots , c^{I+1} \end{array}\right]^T
,
\qquad
\mathbf{f} = \left[\begin{array}{c}-\ui|_{\Gamma_1} ,
-\frac{\partial \ui}{\partial \mathbf{n}}|_{\Gamma_1} , 0 , \cdots , 0 \end{array}\right]^T
~.
\end{equation}
Then all of the coupled BIEs and functional equations
\eqref{bie_interface1}-\eqref{bie_interface6}
can be compactly grouped into the matrix-type notation,
\begin{equation}
\mathbf{A}\bm{\eta}+\mathbf{B} \mathbf{c} = \mathbf{f}
~, \label{rokhlin_muller_matrix}
\end{equation}
where $\mathbf{A}$ is a $I$-by-$I$ matrix, each of whose entries $\mbf{A}_{i,j}$ is
a $2\times2$ block
of operators which maps $\eta_j$ to a pair of functions
(i.e. values then normal derivatives) on $\Gamma_i$.
Every block of $\mbf{A}$ is zero apart from the following tridiagonal entries,
\begin{eqnarray}
\mbf{A}_{i,i} &= &\left[\begin{array}{cc}-\mathbf{I}+(\tilde{D}^i_{\Gamma_i,\Gamma_i}-\tilde{D}^{i+1}_{\Gamma_i,\Gamma_i}) & (\tilde{S}^i_{\Gamma_i,\Gamma_i}-\tilde{S}^{i+1}_{\Gamma_i,\Gamma_i}) \\
(\tilde{T}^i_{\Gamma_i,\Gamma_i}-\tilde{T}^{i+1}_{\Gamma_i,\Gamma_i}) & \mathbf{I}+(\tilde{D}^{i,*}_{\Gamma_i,\Gamma_i}-\tilde{D}^{i+1,*}_{\Gamma_i,\Gamma_i})\end{array}\right],\quad i=1,2,\cdots,I,
\nonumber\\
\mbf{A}_{i,i+1} &=& \left[\begin{array}{cc}-\tilde{D}^{i+1}_{\Gamma_i, \Gamma_{i+1}} & - \tilde{S}^{i+1}_{\Gamma_{i}, \Gamma_{i+1}} \\
 -\tilde{T}^{i+1}_{\Gamma_{i}, \Gamma_{i+1}} & - \tilde{D}^{{i+1},*}_{\Gamma_{i}, \Gamma_{i+1}}\end{array}\right], \qquad\quad i=1,2,\cdots,I-1,
\nonumber\\
\mbf{A}_{i,i-1} &=& \left[\begin{array}{cc}\tilde{D}_{\Gamma_{i},\Gamma_{i-1}}^{i} & \tilde{S}_{\Gamma_{i},\Gamma_{i-1}}^{i} \\
\tilde{T}_{\Gamma_{i},\Gamma_{i-1}}^{i} & \tilde{D}_{\Gamma_{i},\Gamma_{i-1}}^{{i},*}\end{array}\right], \qquad\qquad i=2,3,\cdots,I
\label{Ablocks}
~,
\end{eqnarray}
where $\mbf{I}$ is the identity operator.
$\mathbf{B}$ is an $I$-by-$(I+1)$ matrix, each of whose entries $\mbf{B}_{i,j}$ is a
stack of $P$ continous function columns (sometimes called a {\em quasi-matrix})
expressing the
effect of each proxy point strength $\cp{j}$ on the value and normal derivative
functions on $\Gamma_i$.
The only nonzero blocks of $\mbf{B}$ are
\be
\mbf{B}_{i,i} =  \left[\begin{array}{ccc} \phi^i_1|_{\Gamma_i},&\cdots,&\phi^i_P|_{\Gamma_i}\\
\bigl.\frac{\partial\phi^i_1}{\partial \nn}\bigr|_{\Gamma_i},&\cdots,
&\bigl.\frac{\partial\phi^i_P}{\partial \nn}\bigr|_{\Gamma_i}
\end{array}\right],
\;
\mbf{B}_{i,i+1} =  \left[\begin{array}{ccc} -\phi^{i+1}_1|_{\Gamma_i},&\cdots,&-\phi^{i+1}_P|_{\Gamma_i}\\
-\bigl.\frac{\partial\phi^{i+1}_1}{\partial \nn}\bigr|_{\Gamma_i},&\cdots,
&-\bigl.\frac{\partial\phi^{i+1}_P}{\partial \nn}\bigr|_{\Gamma_i}
\end{array}\right], \; i=1,2,\cdots, I
\label{Bblocks}
\ee
The term $\mbf{A}\bm{\eta}$ in \eqref{rokhlin_muller_matrix} is precisely
(barring the summation over neighbors) the
M\"uller--Rokhlin formulation \cite{muller,rokh83} for multiple material interfaces.
This is of Fredholm second kind since
the off-diagonal blocks in \eqref{Ablocks} have continuous kernels,
and cancellation of the leading singularities occurs in the pairs
in parentheses in \eqref{Ablocks},
leaving the diagonal operators at most weakly singular, hence compact.

\subsection{Imposing the quasi-periodicity conditions} \label{qc}

Quasi-periodicity \eqref{QP}
will be enforced in each layer by matching both values and normal
derivatives between the left $L_i$ and right $R_i = L_i+\mbf{d}$ walls.
Since the PDE is second-order, matching two functions (values and normal derivatives)
is sufficient Cauchy data to guarantee extension as a quasi-periodic solution.

We evaluate the first layer representation \eqref{firstlayer_sol} on the walls,
and exploit the following simplification due to translational symmetry
(as in \cite{qpsc, qplp}) which cancels
six terms (three from each near-field sum) down to two,
\bea
&&\hspace{-.2in}\al^{-1}u_1|_{R_1} - u_1|_{L_1} =
\al^{-1}\left(\tilde{D}^1_{R_1, \Gamma_1} \tau_1+ \tilde{S}^1_{R_1,\Gamma_1} \sigma_1 +
\sump \cphi{1}|_{R_1}\right)
\;-\;\left(\tilde{D}^1_{L_1, \Gamma_1} \tau_1+ \tilde{S}^1_{L_1,\Gamma_1} \sigma_1 +
\sump \cphi{1}|_{L_1}\right)
\nonumber\\
&&=\left(\al^{-2} D^1_{R_1+\dd,\Gamma_1}- \al D^1_{L_1-\dd,\Gamma_1}\right) \tau_1 +
\left(\al^{-2} S^1_{R_1+\dd,\Gamma_1}- \al S^1_{L_1-\dd,\Gamma_1}\right) \sigma_1 +
\sump \left( \al^{-1}\phi^1_p|_{R_1} - \phi^1_p|_{L_1}\right) c^1_p
\eea
For quasi-periodicity we wish this function to vanish, so we make it the first
operator block row of a homogeneous linear system.
Doing the same for the normal derivatives on the $L_1$ and $R_1$ walls,
and then for similar conditions for all other layers $i=2,\ldots,I+1$,
gives  equations that can be written compactly with a matrix notation as follows:
\begin{equation}
\mathbf{C}\bm{\eta}+\mathbf{Q}\mathbf{c} = \mathbf{0}~,
\end{equation}
where $\mathbf{C}$ is an $(I+1)$-by-$(I+1)$ matrix, each entry of which is a $2\times2$
block of operators mapping interface densities to wall values and normal derivatives.
Every block of $\mbf{C}$ is zero apart from the bidiagonal blocks,
\bea
\mbf{C}_{i,i} &=& \left[\begin{array}{cc}
             \alpha^{-2} D^i_{R_i+\mathbf{d}, \Gamma_i}- \alpha D^i_{L_i-\mathbf{d}, \Gamma_i} \;&
             \alpha^{-2} S^i_{R_i+\mathbf{d},\Gamma_i}- \alpha S^i_{L_i-\mathbf{d}, \Gamma_i} \\
             \alpha^{-2} T^i_{R_i+\mathbf{d}, \Gamma_i}- \alpha T^i_{L_i-\mathbf{d}, \Gamma_i}\; & 
             \alpha^{-2} D^{i,*}_{R_i+\mathbf{d},\Gamma_i}- \alpha D^{i,*}_{L_i-\mathbf{d}, \Gamma_i}\end{array}\right]
\label{Cii}
\\
\mbf{C}_{i,i-1} &=& \left[\begin{array}{cc}
             \alpha^{-2} D^i_{R_i+\mathbf{d}, \Gamma_{i-1}}- \alpha D^i_{L_i-\mathbf{d}, \Gamma_{i-1}} \;&
             \alpha^{-2} S^i_{R_i+\mathbf{d},\Gamma_{i-1}}- \alpha S^i_{L_i-\mathbf{d}, \Gamma_{i-1}} \\
             \alpha^{-2} T^i_{R_i+\mathbf{d}, \Gamma_{i-1}}- \alpha T^i_{L_i-\mathbf{d}, \Gamma_{i-1}}\; & 
             \alpha^{-2} D^{i,*}_{R_i+\mathbf{d},\Gamma_{i-1}}- \alpha D^{i,*}_{L_i-\mathbf{d}, \Gamma_{i-1}}\end{array}\right]
\eea
for $i=1,2,\cdots,I+1~$ and $i=2,3, \cdots, I+1,$ respectively.
 $\mathbf{Q}$ is an $(I+1)$-by-$(I+1)$ matrix, each entry of which is
a stack of $P$ function columns (as with $\mbf{B}_{i,j}$),
but only the diagonal entries are nonzero,
\be
\mbf{Q}_{i,i} \;=:\; \mbf{Q}_i \;=\; \left[\begin{array}{ccc}
 \al^{-1}\phi^i_1|_{R_i} - \phi^i_1|_{L_i} \;, &\cdots, &
 \al^{-1}\phi^i_P|_{R_i} - \phi^i_P|_{L_i}
\\
 \al^{-1}\bigl.\frac{\partial\phi^i_1}{\partial\nn}\bigr|_{R_i} - \bigl.\frac{\partial\phi^i_1}{\partial\nn}\bigr|_{L_i} \;,&\cdots,&
 \al^{-1}\bigl.\frac{\partial\phi^i_P}{\partial\nn}\bigr|_{R_i} - \bigl.\frac{\partial\phi^i_P}{\partial\nn}\bigr|_{L_i}
\end{array}\right]
\mbox{ for } i=1,2,\cdots,I+1~.
\label{Qiop}
\ee

\subsection{Imposing the radiation conditions} \label{rc}

First we enforce the
upward radiation condition \eqref{RB1}
at the artificial interface $U$ (with upward-pointing normal),
%
%
substituting the layer-1 representation
\eqref{firstlayer_sol} to get,
\be
\tilde{D}^1_{U, \Gamma_1} \tau_1+ \tilde{S}^1_{U,\Gamma_1} \sigma_1 +
\sump c^1_p \phi^1_p|_U
- \sum_{n \in \mathbb{Z}} a^U_n e^{i\kappa_n x} \;=\; 0~.
\ee
Matching values at $U$ is not enough: we also need to match normal ($y$) derivatives,
to ensure that the second-order PDE solution continues smoothly through $U$, thus,
\be
\tilde{T}^1_{U, \Gamma_1} \tau_1+ \tilde{D}^{1,\ast}_{U,\Gamma_1} \sigma_1 +
\sump c^1_p \biggl.\frac{\partial\phi^1_p}{\partial\nn}\biggr|_U
- \sum_{n \in \mathbb{Z}} a^U_n ik^U_n e^{i\kappa_n x} \;=\; 0~.
\ee
We will truncate the Rayleigh--Bloch expansion to $2\nrb+1$ terms,
from $n=-\nrb$ to $\nrb$,
since it is exponentially convergent once $|\kappa_n|$ exceeds $k_1$ (in the upper
layer) and $k_{I+1}$ (lower layer).
We also apply the downward radiation condition \eqref{RB3} at $D$ to the representation
\eqref{lastlayer_sol}, giving a second set of homogeneous linear conditions.
We choose the normals of $U$ and $D$ both to point in the upward sense.
As with $\bm{\eta}$ and $\mbf{c}$, we stack all coefficients into a single vector,
\be
\mathbf{a} =
\left[ \mathbf{a}^U, \mathbf{a}^D\right]^T =
\left[a^U_{-\nrb},\cdots, a^U_\nrb,  a^D_{-\nrb},\cdots, a^D_\nrb\right]^T~.
\ee
The resulting conditions can again be written in a simple matrix form:
\begin{equation}
\mathbf{Z} \bm{\eta}+\mathbf{V}\mathbf{c}+\mathbf{W} \mathbf{a} = \mathbf{0}~,
\end{equation}
where
\begin{equation}
\mathbf{Z} = \left[\begin{array}{cccc}Z_{U} & 0 & \cdots & 0\\
                                        0 & \cdots & 0 & Z_D\end{array}\right],~
\mathbf{V} = \left[\begin{array}{cccc}V_{U} & 0 & \cdots & 0\\
                                        0 & \cdots & 0 & V_D\end{array}\right],~
\mathbf{W} = \left[\begin{array}{cc}W_U & 0 \\0 & W_D\end{array}\right],
\end{equation}
in which
\bea
\mbf{Z}_U &=& \left[\begin{array}{cc}\tilde{D}^1_{U,\Gamma_1} & \tilde{S}^1_{U,\Gamma_1}\\
					  \tilde{T}^1_{U,\Gamma_1} & \tilde{D}^{1,*}_{U,\Gamma_1}
\end{array}\right],~
\mbf{Z}_D = \left[\begin{array}{cc}\tilde{D}^{I+1}_{D,\Gamma_I} & \tilde{S}^{I+1}_{D,\Gamma_I}\\
					  \tilde{T}^{I+1}_{D,\Gamma_I} & \tilde{D}^{{I+1},*}_{D,\Gamma_I}
\end{array}\right],
\\
\mbf{V}_U &=& \left[\begin{array}{ccc}\phi^1_1|_U
~,&\cdots,&
\phi^1_P|_U \\
\bigl.\frac{\partial\phi^1_1}{\partial\nn}\bigr|_U
~,&\cdots,&
\bigl.\frac{\partial\phi^1_P}{\partial\nn}\bigr|_U
\end{array}\right],~
\mbf{V}_D = \left[\begin{array}{ccc}\phi^{I+1}_1|_D
~,&\cdots,&
\phi^{I+1}_P|_D \\
\bigl.\frac{\partial\phi^{I+1}_1}{\partial\nn}\bigr|_D
~,&\cdots,&
\bigl.\frac{\partial\phi^{I+1}_P}{\partial\nn}\bigr|_D
\end{array}\right],
\label{VUD}
\\
\mbf{W}_U \!\!&= &\!\!\!
\left[\begin{array}{ccc} -e^{i\kappa_{-\nrb}x}|_U\;,  & \cdots ,&  -e^{i\kappa_{\nrb}x}|_U   \\ -ik^U_{-\nrb}e^{i\kappa_{-\nrb}x}|_U\;, & \cdots,& -ik^U_{\nrb}e^{i\kappa_{\nrb}x}|_U\end{array}\right],
\mbf{W}_D = \left[\begin{array}{ccc} -e^{i\kappa_{-\nrb}x}|_D \;, & \cdots, &  -e^{i\kappa_{\nrb}x}|_D   \\ ik^D_{-\nrb}e^{i\kappa_{-\nrb}x}|_D \;,& \cdots,& ik^D_\nrb e^{i\kappa_{\nrb}x}|_D\end{array}\right]~.
\label{WUD}
\eea
To clarify, in $\mbf{W}_U$ and $\mbf{W}_D$, the $2\nrb+1$
columns are pairs of Fourier functions
evaluated over $x\in(-d/2,d/2)$, the $x$-coordinate extent of the lines $U$ and $D$.

\subsection{Discretization of functions and operators}
Finally, combining the linear conditions from the previous three subsections,
we have the coupled BIE and functional equations,
 \begin{equation}
 \left[\begin{array}{ccc}\mathbf{A} & \mathbf{B} & \mathbf{0} \\ \mathbf{C} & \mathbf{Q} & \mathbf{0} \\ \mathbf{Z} & \mathbf{V} & \mathbf{W}\end{array}\right] \left[\begin{array}{c}\bm{\eta} \\ \mathbf{c} \\ \mathbf{a} \end{array}\right] \;=\; \left[\begin{array}{c}\mathbf{f} \\ \mathbf{0} \\ \mathbf{0}\end{array}\right]~,
\label{system_bie}
 \end{equation}
Recall that $\bm{\eta}$ contains unknown density functions, while $\mbf{c}$ and $\mbf{a}$
are discrete coefficient vectors. On the right hand side, $\mbf{f}$ involves functions
(from the incident wave), and each $\mbf{0}$ is a stack of zero functions.
Thus $\mbf{A}$, $\mbf{C}$, and $\mbf{Z}$ are blocks of operators, while the
other six matrix blocks involve quasi-matrices (stacks of function columns).

For numerical computation the continuous variables must be discretized, turning
each function into a discrete set of values, and each operator block into a matrix.
This is simple for the functional conditions in the
2nd and 3rd block rows of \eqref{system_bie}: we just sample at a discrete
set of collocation points. For the 2nd block row,
we use $M_w$ nodes $\{\xx^i_m\}_{m=1}^{M_w}$
of a Gauss--Legendre quadrature living on the left wall $L_i$ for the $i$th layer.
See \fref{domain}(b).
Thus each diagonal block of $\mbf{Q}$ \eqref{Qiop} is replaced by a
$2M_w$-by-$2P$ matrix $Q_i$ with elements
\be
(Q_i)_{mp} \;=\; \left\{\begin{array}{ll}
 \al^{-1}\phi^i_p(\xx^i_m+\dd) - \phi^i_p(\xx^i_m)~, &
m = 1,\ldots,M_w, \quad p=1,\ldots, P
\\
 \al^{-1}\frac{\partial\phi^i_p}{\partial\nn}(\xx^i_{m-M_w}+\dd) - \frac{\partial\phi^i_p}{\partial\nn}(\xx^i_{m-M_w})~, &
m = M_w+1,\ldots,2M_w, \quad  p=1,\ldots, P
\end{array}\right.
\label{Qi}
\ee
Similarly for the 3rd row, we use $M$ equally-spaced (trapezoid rule)
nodes $\{\xx^U_m\}_{m=1}^M$
on $U$, and $\{\xx^D_m\}_{m=1}^M$ on $D$.
The trapezoid rule is appropriate here since functions will be periodic.
Inserting these nodes, the formulae for the matrices $V_U$, $V_D$, $W_U$, and $W_D$
discretizing \eqref{VUD}-\eqref{WUD} are similar to \eqref{Qi} above.

To discretize the remaining blocks $\mbf{A}$, $\mbf{B}$, $\mbf{C}$ and $\mbf{Z}$,
we need to fix a set of quadrature nodes $\{\zz^i_j\}_{j=1}^{N_i}$ on each
interface $\Gamma_i$.
These nodes have associated weights $\{w^i_j\}_{j=1}^{N_i}$, such that
for any smooth $d$-periodic function $f$ on $\Gamma_i$,
the quadrature rule
\be
\int_{\Gamma_i} f(\rr) ds_\rr \; \approx \; \sum_{j=1}^{N_i} f(\zz^i_j) w^i_j
\label{Giquad}
\ee
holds to high accuracy.
To choose these nodes and weights, 
we first consider the case of $\Gamma_i$ a smooth interface
(e.g.\ $\Gamma_1$ in \fref{domain}(b)).
Let one period of the interface be
parametrized by a vector function $Z(s)$ for $0\le s < 2\pi$.
By changing variable, the
periodic trapezoid rule $s_j = 2\pi (j-1/2)/N_i$, $j=1,\ldots,N_i$ in parameter $s$
gives a quadrature rule $\zz^i_j = Z(s_j)$ and $w^i_j = (2\pi/N_i)|Z'(s_j)|$.
Then for ($C^\infty$) smooth 
$d$-periodic integrands on $\Gamma_i$ the error in \eqref{Giquad}
will be superalgebraically convergent \cite[(2.9.16)]{davisrabin}.
On the other hand, if $\Gamma_i$ has corners, it breaks into one or more ``segments''
(e.g.\ $\Gamma_2$ in \fref{domain}(b)). These need not be straight lines,
merely smooth.
Say a segment is again parametrized by a function $Z(s)$ for $0\le s < 2\pi$.
Then we {\em reparametrize} it via $Z(w(s))$ using the corner
grading function suggested by Kress \cite[(6.9)]{kress91},
$$
w(s) = 2\pi \frac{v(s)^q}{v(s)^q + v(2\pi-s)^q}~,
\qquad \mbox{ where } \quad v(s) = \left(\frac{1}{q}-\frac{1}{2}\right)
\left(\frac{\pi-s}{\pi}\right)^3 + \frac{1}{q}\frac{s-\pi}{\pi} + \frac{1}{2}
~,
\qquad
 0\le s< 2\pi~,
$$
where $q$ controls the grading at endpoints. Higher $q$ will cause more nodes
to be close to the endpoints; typically we choose $q=6$ or higher.
Let $n_{i,l}$ be the number of nodes used for the $l$th segment of $\Gamma_i$.
Then a separate trapezoid rule $s_j = 2\pi j/n_{i,l}$, $j=1,\ldots,n_{i,l}$
is used on each segment, making $N_i = \sum_{l} n_{i,l}$ nodes in total.
The formula for the nodes and weights%
\footnote{We implement this in the \mpspack\ command {\tt segment(...,'pc')}.}
are the same as in the smooth case,
with the proviso that $Z$ is replaced by the composed function $Z \circ w$.
The grading function, since its derivative vanishes to high order at the
endpoints, insures that the trapezoid rule achieves high order accuracy
(typically order $q$).
We find this more efficient for up to 10-digit accuracy than
dyadically-refined panel quadratures \cite{triplejuncQPFDS}.

With interface quadratures defined, blocks $\mbf{B}$, $\mbf{C}$, and $\mbf{Z}$
are easy to discretize by restriction of the continuous variable to the set of nodes.
For example, each block $\mbf{B}_{i,i}$ in \eqref{Bblocks}, describing the
interaction of the $i$th proxy basis with $\Gamma_i$,
is replaced by a $2N_i$-by-$P$ matrix $B_{i,i}$ with elements
\be
(B_{i,i})_{jp} \;=\; \left\{\begin{array}{ll}
 \phi^i_p(\zz^i_j)~, &
j = 1,\ldots,N_i, \quad p=1,\ldots, P
\\
 \frac{\partial\phi^i_p}{\partial\nn}(\zz^i_{j-N_i})~, &
j = N_i+1,\ldots,2N_i, \quad  p=1,\ldots, P
\end{array}\right.
\label{Bblockmat}
\ee
The matrix for $\mbf{B}_{i,i+1}$ is similar.
Operator blocks $\mbf{C}$ and $\mbf{Z}$ involve boundary integral representations
over interfaces: each integral is replaced by a sum according to \eqref{Giquad}.
For example, the $M_w$-by-$N_i$ matrix discretizing
the upper-right block of $\mbf{C}_{i,i}$ in \eqref{Cii}
has elements
$\bigl(\al^{-2} G^i(\xx^i_m+\dd,\zz^i_j) - \al G^i(\xx^i_m-\dd,\zz^i_j) \bigr)w^i_j$,
for $m=1,\ldots,M_w$ and $j=1,\ldots,N_i$.
Other blocks of $C$ and $Z$ are discretized similarly; for the reader's
sanity we refrain from giving all formulae.

Finally we discretize $\mbf{A}$.
The operators in the blocks $\mbf{A}_{i,j}$ for $i\neq j$,
and for the neighboring terms $l=\pm 1$ in the local sums
\eqref{SDopsum}--\eqref{DtTopsum} even when $i=j$, 
involve only interactions between differing interfaces, and thus
may be replaced simply by substitution of the native quadrature rule
\eqref{Giquad} for the sources, and evaluation at discrete target nodes, as above.
This method of discretizing an integral operator is called Nystr\"om's method
\cite[Sec.~12.2]{LIE}.
This leaves only the {\em self-interaction} terms $l=0$ in $\mbf{A}_{i,i}$,
which involve operators that are logarithmically singular.
To achieve high-order accuracy for these operators,
we use the standard ``plain'' Nystr\"om matrix entries
of the form $A(\zz^i_m,\zz^i_j) w^i_j$, for $m,j=1,\ldots,N_i$
(here $A$ symbolizes a generic operator block)
for entries far from the diagonal.
Near-diagonal entries are adjusted by local Lagrange interpolation of the
smooth density from the existing periodic trapezoid nodes onto a
set of auxiliary nodes special to the logarithmic singularity,
due to Alpert \cite{alpert}; see \cite[Sec.~4]{hao} for the full recipe.
We use 30 auxiliary nodes per target node, which achieves high-order
convergence with error
${\mathcal O}(N_i^{-16} \log N_i)$.
With this Nystr\"om matrix, the
non-zero right-hand side terms in $\mbf{f}$ become the samples at the first interface nodes $\zz^1_j$,
and all the unknown vectors $\eta_i$ become samples
of the densities $\tau_i$ and $\sigma_i$ at the nodes of all interfaces.

The size of the resulting matrix $\mbf{A}$ is $\Nden$-by-$\Nden$, where the total number
of density unknowns is
\be
\Nden := 2 \sum_{i=1}^I N_i ~,
\label{Nden}
\ee
the factor of two coming from the two types of layer potential per interface.

\begin{rmk} 
The Alpert correction to the periodic trapezoid rule \cite{alpert}
is designed for closed curves, where the kernel and densities are periodic.
However, our interfaces $\Gamma_i$ do not close on each other, moreover
the solution density $\eta_i$ is {\em quasi-periodic} with Bloch phase $\al$.
This means that,
when $\Gamma_i$ is smooth and hence has a single periodic quadrature rule,
we must modify the Alpert correction entries in the northeast
and southwest corners of the matrix, to account for the continuation of the
interface into the next unit cell, and the phase factors $\al$ and $\al^{-1}$.%
\footnote{This periodic-segment adjustment of the Alpert correction is
in the \mpspack\ code {\tt \@quadr.alpertizeselfmatrix}.}

Kress' periodic logarithmic correction would also be a slightly more
accurate option
\cite{kress91} \cite[Sec.~6]{hao}; however, we found that it was less convenient
to adjust this scheme for quasi-periodic densities and open interface segments.
(See Meier et al.\ \cite{CWnystrom}
for an example of this; extra cut-off functions are
required.)
\label{phasewrap}
\end{rmk}


We will see that, due to the high-order convergence,
the numbers of collocation nodes $N_i$, $n_{i,l}$, $M_w$, and $M$
can be small, of order a hundred, even for 10-digit accuracies.
Note that for smooth interfaces
the periodic trapezoid rule is in fact inaccurate for the interactions
from neighboring interfaces, e.g.\ terms like $S^i_{\Gamma_i\pm\dd,\Gamma_i}$,
due to ``dangling'' ends of these interfaces, but that this is handled by
the periodizing scheme, retaining high-order convergence.


To summarize, the linear algebraic system which discretizes
\eqref{system_bie} has identical structure to \eqref{system_bie},
but with all of its blocks matrices as constructed above.
We notate these blocks using standard (non-bold) font.
It will now be rearranged in order to solve it in a fast direct fashion,
exploiting its tridiagonal structure.

\section{Rearrangement of equations, Schur complement, and fast solver}
\label{solve}

The discretized BIEs of the previous section require a few hundred to a couple
of thousand
unknowns per layer, to represent typical geometries as shown in \fref{domain}
to accuracies of around 10 digits.
The total number of unknowns includes the densities, proxy strengths, and scattered
amplitudes, namely
$$
\Ntot = \Nden + IP + 2(2K+1)
~.
$$
We will find $P \sim 10^2$, thus there is an order of magnitude less
proxy unknowns than density unknowns.
In any case, a device with dozens or more layers leads to linear systems that are
too large for direct $\bigO(\Ntot^3)$ inversion or solution.
On the other hand, due to the proxy point (MFS)
representation, the full linear system is exponentially ill-conditioned
\cite{mfs}, so iterative solution of the full system is impossible.
Therefore, for robustness, we describe in this section a direct solution technique,
that will be ``fast'' (optimally scaling in $I$ the number of layers),
by exploiting the algebraic structure.
This also will make it easier to solve for the technologically-important
case of multiple incident angles $\ti$
at the same frequency $\omega$.

\subsection{Rearrangement}\label{rearrangement}
The first step is to rearrange the unknowns, in a form amenable to elimination
of the proxy and scattering amplitude unknowns.
We reorder our vector of all unknowns to be
\begin{equation}
\mathbf{x}= \left[ \eta, c_1, \mathbf{a}^U, c_2, c_3, \cdots c_{I-1}, c_I, c_{I+1}, \mathbf{a}^D \right]^T = [\mathbf{\eta}, \mathbf{x}_1, \mathbf{x}_2, \cdots, \mathbf{x}_{I+1}]^T
\end{equation}
where
\be
\mathbf{x}_1 := \left[c_1, \mathbf{a}^U \right]^T, \qquad
\mathbf{x}_i := c_i, \quad i=2,3,\cdots,I~, \qquad
\mathbf{x}_{I+1} := \left[c_{I+1}, \mathbf{a}^D \right]^T.
\ee
Now similarly rearranging the (discretized) blocks of 
\eqref{system_bie} we get the full linear system,
 \begin{equation}
 \left[\begin{array}{cccccccccccc}
   & \multicolumn{5}{c}{\multirow{5}*{  $A$   }}
   & B'_{1,1}   & B_{1,2} & 0 & \cdots &0 \\
  &&&&&& 0 & B_{2,2}  & B_{2,3}& \cdots  &0\\
  &&&&&   & 0 & 0 & B_{3,3}& \cdots&0 \\
  &&&&&&\vdots&\vdots &\vdots&\vdots &\vdots\\
  &&&&& &0& 0& 0 & \cdots  &B_{I-1,I+1}\\
  &&&&& &0& 0& 0&  \cdots  & B'_{I,I+1}\\
 C'_{1,1} & 0 & 0 &\cdots  &0 &0 & Q'_1& 0& 0& \cdots &0\\
 C_{2,1} & C_{2,2} & 0 & \cdots &0 &0 &0& Q_2& 0& \cdots &0\\ 
0 & C_{3,2} & C_{3,3} &\cdots &0  &0 &0& 0& Q_3& \cdots &0\\ 
\vdots & \vdots & \vdots &\vdots  &\vdots &\vdots  &\vdots&\vdots&\vdots &\vdots&\vdots\\ 
0 & 0 &  0 &\cdots&C_{I,I-1} &C_{I,I} &0 &  0& 0& \cdots &0\\ 
0 & 0 & 0 &\cdots&0 & C'_{I+1, I} &0 & 0&  0& \cdots &Q'_{I+1}\\ 
 \end{array}\right] \mathbf{x} = \left[\begin{array}{c}\mathbf{f} \\ \mathbf{0} \\\mathbf{0} \\ \vdots \\\mathbf{0} \\\mathbf{0} \\\mathbf{0} \\\mathbf{0} \\\mathbf{0} \\\vdots \\\mathbf{0} \\\mathbf{0}\end{array}\right].
\label{linsys}
 \end{equation}
Because the amplitudes $\mathbf{a}$ were separated into $\mathbf{a}^U$ and $\mathbf{a}^D$
then merged into $c_1$ and $c_{I+1}$ to make $\mathbf{x}_1$ and $\mathbf{x}_{I+1}$,
respectively, the first and last blocks of $B$, $C$ and $Q$ matrix had to be expanded
to
\begin{eqnarray}
 B'_{1,1} &=& \left[\begin{array}{cc}B_{1,1} & \mathbf{0}\end{array}\right],
\qquad B'_{I,I+1} = \left[\begin{array}{cc}B_{I,I+1} & \mathbf{0}\end{array}\right],
\nonumber
\\
C'_{1,1}  &=& \left[\begin{array}{c}C_{1,1} \\Z_U\end{array}\right],
\qquad
C'_{I+1,I} = \left[\begin{array}{c}C_{I+1,I} \\Z_U\end{array}\right],
\nonumber
\\
Q'_{1}  &=& \left[\begin{array}{cc}Q_1 & \mathbf{0} \\ V_U & W_U\end{array}\right],
\qquad
Q'_{I+1}  = \left[\begin{array}{cc}Q_{I+1} & \mathbf{0} \\ V_D & W_D\end{array}\right].
\nonumber
\end{eqnarray}

\subsection{Schur complements}
\label{schur}
The rearrangement in subsection \ref{rearrangement} enables us to use $I+1$ independent
Schur complements to eliminate all the unknown vectors $\mathbf{x}_i$;
this corresponds to periodizing all of the Green's functions.
For example, consider the equation in the first row,
$$
A_{1,1} \eta_1+A_{1,2} \eta_2+B'_{1,1} \mathbf{x}_1+B_{1,2} \mathbf{x}_2 = \mathbf{f}
~.
$$
The first two equation rows starting the $\mathbf{C}$ block are
$$
C'_{1,1}\eta_1+Q'_{1}\mathbf{x}_1 = \mathbf{0}~, \qquad
C_{2,1}\eta_1+C_{2,2}\eta_2+Q_2\mathbf{x_2} = 0~.
$$
With the first of these $\mathbf{x}_1$ can be eliminated,
and with the second $\mbf{x}_2$ can, giving
\be
(A_{1,1}-B'_{1,1}Q^{'\dag}_1 C'_{1,1}  - B_{1,2}Q^{\dag}_2 C_{2,1})\eta_1+(A_{1,2}- B_{1,2}Q^{\dag}_2 C_{2,2})\eta_2 = \mathbf{f}
~,
\label{elim}
\ee
where $Q^{'\dag}_i$ denotes the {\em pseudo-inverse} of the rectangular matrix
$Q'_i$, for $i=1,\ldots,I+1$.
Filling the matrix $Q^{'\dag}_i$ and then using matrix-matrix multiplication
with the $C$ blocks to fill matrices in \eqref{elim} would lose
accuracy, due to round-off error with the exponentially large matrix entries of the pseudo-inverses.
To retain full accuracy, we do not in fact ever evaluate the pseudo-inverses.
Rather, for product matrices such as $X = Q^{'\dag}_1 C'_{1,1}$ appearing above,
we solve the (ill-conditioned) linear system
$ Q'_1 X = C'_{1,1}$ with a standard backward-stable direct dense solver.%
\footnote{For all such dense solves we
use the ``backslash'' or {\tt mldivide} command in MATLAB.}

Similar computations eliminate $\mathbf{x}_1$, $\mathbf{x}_2$, and $\mathbf{x}_3$ from the second equation to give
$$
(A_{2,1}-B_{2,2}Q^{\dag}_2 C_{2,1})\eta_1+(A_{2,2} -B_{2,2}Q^{\dag}_2 C_{2,2} - B_{2,3}Q^{\dag}_3 C_{3,2} )\eta_2 +(A_{2,3}-  B_{2,3}Q^{\dag}_3 C_{3,3} )\eta_3 = \mathbf{0}
~.
$$
By repeating the same computation for all the equations, all the $\mathbf{x}_j$ are eliminated and a block tridiagonal system for $\eta$ is obtained as
\begin{equation}
 \left[\begin{array}{cccccccc}
 A'_{1,1} & A'_{1,2} & 0 & 0 & \cdots &0 & 0 & 0\\
 A'_{2,1} & A'_{2,2} & A'_{2,3} & 0 & \cdots &0 & 0 & 0\\
 0 & A'_{3,2} & A'_{3,3} & A'_{3,4} & \cdots &0 & 0 & 0\\
 \vdots & \vdots & \vdots & \vdots & \vdots & \vdots & \vdots & \vdots\\
 0 & 0 & 0 & 0 & \cdots & A'_{I-1,I-2} & A'_{I-1,I-1} & A'_{I-1,I}\\
 0 & 0 & 0 & 0 & \cdots & 0 & A'_{I,I-1} & A'_{I,I}
 \end{array}\right]\left[\begin{array}{c}\eta_1 \\ \eta_2 \\ \vdots \\ \eta_{I-1} \\ \eta_I\end{array}\right] = \left[\begin{array}{c}\mathbf{f} \\ \mathbf{0} \\\mathbf{0} \\\mathbf{0} \\\mathbf{0}\end{array}\right],
\label{Aschur}
 \end{equation}
 with new interaction matrices for the density unknowns,
 \begin{eqnarray}
 A'_{1,1} &=& A_{1,1}-B'_{1,1}Q_1^{' \dag} C'_{1,1} - B_{1,2}Q_{2}^{\dag}C_{2,1}~,
\label{Schur1}\\
A'_{i,i} &=& A_{i,i}-B_{i,i}Q_i^{\dag} C_{i,i} - B_{i,i+1}Q_{i+1}^{\dag}C_{i+1,i}~,\qquad
i=2,3,\cdots,I-1~,\\
A'_{I,I} &=& A_{I,I}-B_{I,I}Q_I^{\dag} C_{I,I} - B'_{I,I+1}Q_{I+1}^{' \dag}C'_{I+1,I}~,\\
A'_{i,i+1} &=& A_{i,i+1} - B_{i,i+1}Q_{i+1}^{\dag} C_{i+1,i+1}~,\qquad\qquad i=1,2,\cdots,I-1~, \\
A'_{i,i-1} &=&  A_{i,i-1}-B_{i,i}Q_{i}^{\dag}C_{i,i-1}~,\qquad\qquad\qquad i=2,3,\cdots,I~.
\label{Schuri}
 \end{eqnarray}

\begin{rmk}
The Schur complements \eqref{Schur1}-\eqref{Schuri}
replacing matrix blocks $A_{i,j}$ by
$A'_{i,j}$ correspond to replacing each layer free-space Green's
function $G^i$ \eqref{Gi} by an arbitrary quasi-periodic Green's function $G'^i$ that
obeys the layer's correct wall boundary conditions $G'^i(\xx+\dd,\yy)=\al G'^i(\xx,\yy)$,
$\xx\in L_i$, $\yy\in \Omega_i$.
Crucially, since upward and downward radiation conditions are never imposed
together in any single layer $i$, this is {\em not} the standard quasi-periodic $G^i_{QP}$ of
\eqref{GQPi},
which diverges at that layer's Wood anomalies.
The $G'^i$ selected by the backward-stable solves are always non-divergent;
intuitively, since
proxy strength vectors of small norm are possible, hence these are
selected.
\end{rmk}

\subsection{Block tridiagonal solve and evaluation of scattered wave}
\label{fastsolver}

The block tridiagonal system \eqref{Aschur}
can be efficiently solved with a block LU decomposition
\cite[Sec.~4.5.1]{golubvanloan}, at a cost of
dense direct inversion of diagonal blocks.
Since they derive from a second-kind integral equation, these diagonal blocks
are all well conditioned, and no significant rounding error occurs when their
inverses are multiplied as matrices.
We write $f_i$ for the block vectors of the right-hand side of \eqref{Aschur}.
The algorithm is initialized by setting $\tilde{f}_1 = f_1$ and $\tilde{A}_{1,1}=A'_{1,1}$,
then the {\em forward sweep}, for $i=2$ to $I$ in order,
\bea
\tilde{A}_{i,i} &=& A'_{i,i} - A'_{i,i-1} (\tilde{A}_{i-1,i-1})^{-1} A'_{i-1,i} ~,
\nonumber \\
\tilde{f}_{i} &=& f_{i} -   (\tilde{A}_{i-1,i-1})^{-1} A'_{i-1,i} \tilde{f}_i
~.
\nonumber
\eea
To save RAM, it is possible to have $\tilde{A}_{i,i}$ overwrite $A'_{i,i}$.
The {\em backward sweep} starts with solving $\tilde{A}_{I,I} \eta_I = \tilde{f}_I$,
then for $i=I-1$ down to 1 in order, solve for $\eta_i$ in
$$
\tilde{A}_{i,i} \eta_i \; =\; \tilde{f}_i - A'_{i,i+1} \eta_{i+1} ~.
$$
If each $N_i \approx N$, for some constant $N$,
the cost is $\bigO(N^3I)$ due to two block inversions per layer.
This is roughly $I^2$ times faster than naive inversion of the whole system.
Note that the matrix filling time is only $\bigO(N^2I)$,
but dominates for the small $N$ in our settings,
due the large prefactor in evaluating special functions and
applying Alpert corrections.

Once all of the density vectors $\eta_i$ are known,
the proxy and scattering amplitude vectors $\mathbf{x}_i$ are easily recovered by
 \begin{eqnarray}
\mathbf{x}_1 &=& -Q'^{\dag}_1 C'_{1,1}\eta_1~,\\
\mathbf{x}_i &=& -Q^{\dag}_i \left[\begin{array}{cc}C_{i,i-1} & C_{i,i}\end{array}\right]  \left[\begin{array}{c}\eta_{i-1} \\ \eta_{i} \end{array}\right]~,\qquad i=2,3,\cdots,I~,\\
\mathbf{x}_{I+1} &=& -Q'^{\dag}_{I+1} C'_{I+1,I}\eta_I~.
 \end{eqnarray}
Here the products of the type $Q_i^{\dag}C_{i,j}$ can be reused from their
prior computation in the Schur complements \eqref{Schur1}-\eqref{Schuri}.

Finally, the scattered wave solution $u_i$ in each layer can be evaluated
from their representations \eqref{firstlayer_sol}-\eqref{lastlayer_sol}
by applying the interface quadrature rules to the single- and double-layer
potentials, using the discrete density vectors $\eta_i$.
The evaluation involves only free-space Green's functions and monopole/dipole
sources, which would be compatible with a standard FMM \cite{fmm2}.
Above $y=y_U$ and below $y=y_D$, the solution is evaluated
from the Rayleigh--Bloch expansions (truncated versions of \eqref{RB1}-\eqref{RB3}),
using the Bragg amplitudes $\mbf{a}^U$ and $\mbf{a}^D$.

\begin{figure}[th] 
   \centering
   \includegraphics[width=4.5in]{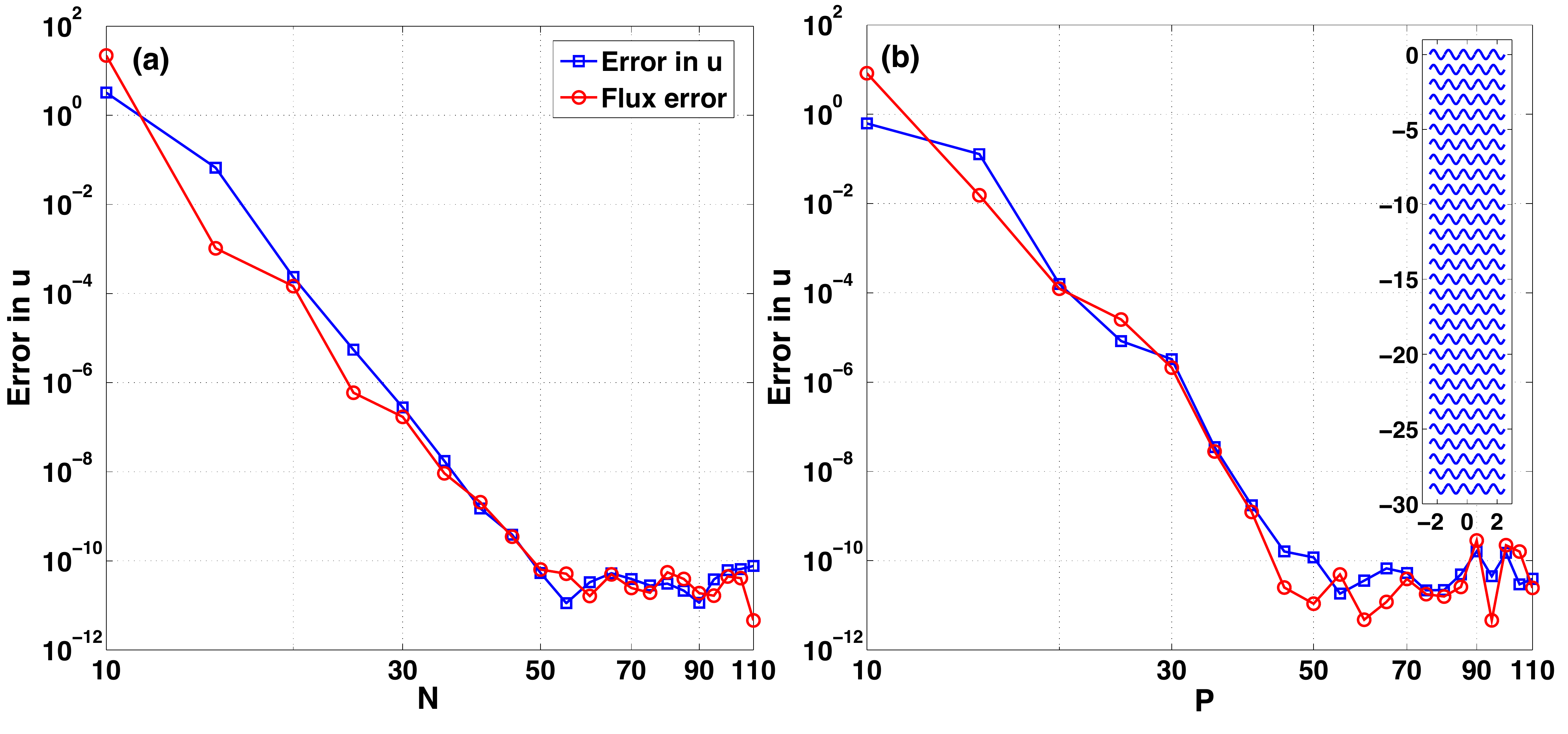} 
   \caption{Convergence of $u(0.15, 0.6)$ and flux error for 30 sine interfaces (see the inset in (b)) with $\omega = 10$ and  random $\varepsilon_i$ between 1 and 2 . {\bf (a)} Convergence in $N$, number of nodes per sine interface (blue square) and flux error (red triangle) while $P = 110$. {\bf (b)} Convergence in $P$ (blue square) and flux error (red triangle), fixing $N = 70$ per sine interface. All other parameters are fixed at $M_w$ = 110, $M$ = 100, $K$ = 10, and $R = 2$. (Color online.)}
   \label{convergence_sine}
\end{figure}
\begin{figure}[th] 
\centering
   \includegraphics[width=6in]{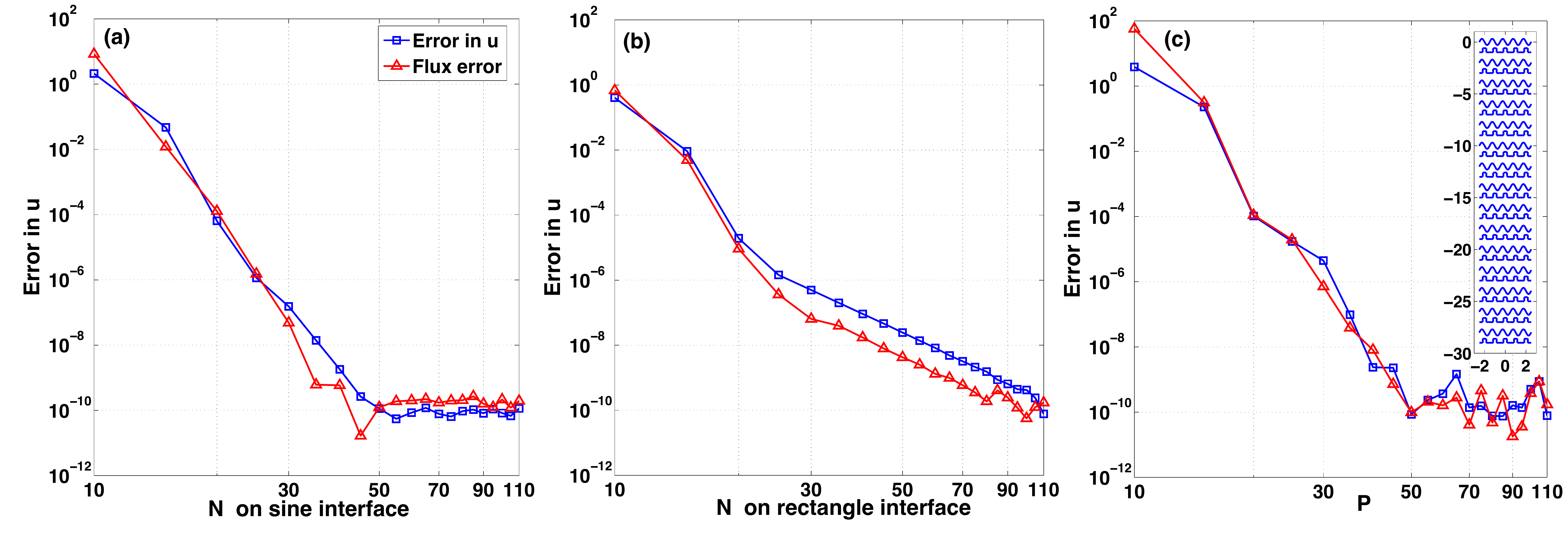} 
   \caption{Convergence of $u(0.15,0.6)$ and flux error for 30 mixed sine and rectangle interfaces (see the inset in (c)) with $\omega = 10$ and random $\varepsilon_i$ between 1 and 2. {\bf (a)} Convergence in $N$ on sine interface (blue square) and flux error (red triangle), while $N = 110$ on each line segment of rectangle interfaces, and $P = 110$. {\bf (b)} Convergence in $N$ on each line segment of rectangle interfaces (blue square) and flux error (red triangle), while $N = 70$ on sine interfaces and $P = 110$. {\bf (c)} Convergence in $P$ (blue square) and flux error (red triangle), while  $N = 70$ on sine, $N=110$ on rectangle interfaces.  All other parameters are fixed at $M_w$ = 110,  $M$ = 120, $K$ = 20, and $R = 2$. (Color online.)}
   \label{convergence_rectangle}
\end{figure}

\subsection{Accelerated sweep over multiple incident angles at one frequency}
\label{sweep}
Modeling many real-world devices requires
characterizing transmission and reflection over a wide range of incident angles $\ti$
at one frequency $\omega$.
In general, changing $\ti$ changes both the operators and the right-hand side in 
\eqref{system_bie}, thus naively an independent solution is required for each
angle, making the task very expensive.
Here we show how to exploit two independent types of structure that speeds
this up by typically an order of magnitude.
  
Firstly, we exploit the fact that
the operators in \eqref{system_bie}, and hence the entire system matrix,
depends on $\ti$ only through $\al$ in \eqref{al}.
Thus we may group together multiple incident angles that share $\al$
(as in \cite{qpfds}), and solve them together at a cost that is
essentially the same as a single angle,
using the precomputed inverse blocks in the tridiagonal solve of \sref{fastsolver}.
Roughly there are $k_1 d/\pi$ such angles, hence this is the speed-up factor.
It grows in proportion to the incident wavenumber.
For example, if $k_1=40$ and $d=1$, the speedup is around a factor of 6.
To set up a sweep of all incident angles, 
we choose a uniform grid in $\ti$ with spacing $2\pi/(nd k_1)$ for some
integer $n$, insuring that only $n$ independent $\al$-value solves are needed.

Secondly, we exploit the fact that filling (rather than solving)
the matrix blocks in \eqref{linsys}
often accounts for the bulk of the solution time.
Consider one of the integral operators used in the BIE matrix (defined in \eqref{SDopsum}),
\begin{equation}
(\Dr^i_{V,W} \tau)(\mathbf{r})= \sum_{l=-1}^1  \alpha^l \int_{W}\frac{\partial G_i}{\partial \mathbf{n}'} (\mathbf{r}, \mathbf{r}'+l \mathbf{d}) \tau(\mathbf{r}') ~ds_{\mathbf{r}'}.
\end{equation}
The integral part is independent of the incident angle or $\alpha$. Therefore
\begin{equation}
\int_{W}\frac{\partial G_i}{\partial \mathbf{n}'} (\mathbf{r}, \mathbf{r}'+l \mathbf{d}) \tau(\mathbf{r}') ~ds_{\mathbf{r}'} \label{precompute}
\end{equation}
can be precomputed for $l\in\{-1,0,1\}$ then used to assemble $(\tilde{D}^i_{V,W} \tau)(\mathbf{r})$ whenever a new $\alpha$ is given.
Exactly the same argument applies to $(\tilde{S}^i_{V,W} \sigma)(\mathbf{r})$,  $(\tilde{T}^i_{V,W} \tau)(\mathbf{r})$, and $(\tilde{D}^{i,*}_{V,W} \sigma)(\mathbf{r})$.
Therefore, the $A$, $B$, $C$, and $Q$ matrices
can be assembled for any $\al$ simply by adding and subtracting
precomputed integral operators; this speeds up the
solution at multiple $\al$ values at the expense of using extra RAM.
%
%

\section{Numerical Results} \label{num}

In all numerical examples, we let the first layer be ``vacuum'' with $\varepsilon_1 = 1$, set $\mu_i=1$ for all layers, and choose periodicity $d = 1$.
All the computations are conducted using MATLAB 2014a running on a workstation with two Intel Xeon E5-2687W processors (total 16 cores) with 256 GB memory and CentOS 6.5.
For filling of the Nystr\"om matrix blocks
we use \mpspack\ \cite{mpspack}, which has an interface to Fortran codes for
Hankel function evaluations by Vladimir Rokhlin. The {\tt parfor} command in MATLAB's Parallel Computing Toolbox
is also used to fill the $A$ and $C$ matrices, using only 12 threads.
%
In the following we present:
numerical tests on convergence, scaling of timing and memory with frequency $\omega$
and the number of layers $I$, solution plots for a 1000-interface case,
and accelerated computation
of transmission and reflection spectra at multiple angles.

\begin{rmk} 
The Bragg coefficients $\{a_n^U\}$ and $\{a_n^D\}$ will be used as an independent
measure of the
accuracy of our numerical scheme based on conservation of the flux (energy)
\cite{linton07,shipmanreview}, namely 
\begin{equation}
\sum_{k^U_n > 0} k^U_n|a^U_n|^2 + \sum_{k^D_n>0}k^D_n|a^D_n|^2 = k_1\cos\ti~.
\end{equation}
This holds when all the material properties $\eps_i$ and $\mu_i$ are real.
Therefore, we will define the relative flux error as
\begin{equation}
\mbox{\rm Flux error }:=\;
\frac{\sum_{k^U_n > 0} k^U_n|a^U_n|^2 + \sum_{k^D_n>0}k^D_n|a^D_n|^2 - k_1\cos\ti}{k_1\cos\ti}~.
\label{fluxerror}
\end{equation}
\label{r:flux}\end{rmk} 

\subsection{Convergence}

\begin{table}[t]
   \centering
\topcaption{\footnotesize{CPU time, memory, and flux error: $\omega = 5$ (period is $0.8\lambda$ in vacuum), $\varepsilon_1 = 1$ and all other $\varepsilon_i$ are random between 1 and 2, $\theta^{inc} = -\pi/5$, $N_i=70$ on sine, $N_i=100\times2$ on triangle, and $N_i=100 \times 5$ on rectangle interfaces, $M_w$ = 120, $M$ = 60, $P = 60$, $K$ = 20, and $R = 2$.}} 
\footnotesize{
   \begin{tabular}{@{} lc c c c c c c l@{}} 
      \toprule
       Number of interfaces 		&  1 & 3  & 10 & 30  & 100 & 300 &1000 \\
      \midrule
       Matrix Filling (sec)    		&   0.518    &  1.860    &  4.200      &  5.600     &   12.384    &    32.332   &  103.331\\
       Schur Complement (sec)    &   0.028     &  0.058    &  0.299     &   0.644     &   2.263      &   6.525    &  21.037\\
       Block Solve (sec)    		&  0.003     &  0.041    & 0.398       &  0.898      &  2.805       &   8.626    &  26.655\\
       Memory (MB)   			&   18         &  41         &   83        &  183         &   608         &   1753       &  5830\\
       Flux Error   			& 4.8e-12   &  3.1e-11 & 2.4e-11    &  4.0e-11   &   2.2e-11   &  1.3e-10    &  9.1e-10 \\
      \bottomrule
   \end{tabular}
   }
   \label{table_omega_5}
\end{table}
\begin{table}[t]
   \centering
\topcaption{\footnotesize{CPU time, memory, and flux error: $\omega = 40$ (period is $6.4\lambda$ in vacuum), $\varepsilon_1 = 1$ and all other $\varepsilon_i$ are random between 1 and 2, $\theta^{inc} = -\pi/5$, $N_i=180$ on sine, $N_i=150\times2$ on triangle, and  $N_i=340\times 5$ on rectangle interfaces, $M_w$ = 120, $M$ = 60, $P = 160$, $K$ = 20, and $R = 2$.} } 
\footnotesize{   \begin{tabular}{@{} l c c c c c c c l@{}} 
      \toprule
       Number of interfaces      &  1 & 3  & 10 &  30 & 100 & 300 &1000 \\
      \midrule
       Matrix Filling (sec)           &  0.813    &   2.739   &   13.295    & 16.312   &  30.960      &  75.915    & 247.597\\
       Schur Complement (sec)  &  0.080    &   0.231   &   1.366      &  3.361    &  10.553      &  33.097     & 108.384\\
       Block Solve (sec)            &  0.018    &   0.118   &    3.879     &  8.637    &   28.180      &  90.187     & 277.933\\
       Memory (MB)   		     & 58          &   112      &    496        &  1069     &   3576         &  10319      & 34347\\
       Flux Error 			     & 8.3e-13  & 5.3e-12  &   8.9e-10  &  3.4e-10 &  3.6e-09      &  5.0e-09   & 4.7e-08\\
      \bottomrule
   \end{tabular}
   }
   \label{table_omega_40}
\end{table}

We study the convergence of the scattered field at a fixed point $(0.15,0.6)$ in the first layer, relative to its converged value, and convergence
of flux errors defined in \eqref{fluxerror}.
The frequency $\omega=10$ corresponds to a period of about 1.6 wavelengths
in the first layer, and larger numbers of wavelengths in the other layers.
Thirty sine interfaces with random $\varepsilon_i$ chosen between 1 and 2 are considered in \fref{convergence_sine}. The actual geometry is depicted in the inset of \fref{convergence_sine}(b). All the interfaces have a periodic trapezoid rule with the same number of quadrature points $N_i = N$, $i=1,2,\cdots 30$.  \fref{convergence_sine}(a) shows convergence of the scattered field error (blue square) and flux error (red circle) as a function of $N$:
the expected 16th-order rate is observed, and the fact that
the two types of error are essentially equivalent.
The same convergence test is conducted as a function of the number of proxy points $P$ per layer, in Fig.~\ref{convergence_sine}(b).
This time super-algebraic convergence is observed.
The fact that 10-digit accuracy results with only $N=P=50$ is a testament
to the extremely rapid convergence of the method.

\begin{figure}[h!] 
   \centering
   \includegraphics[width=5in]{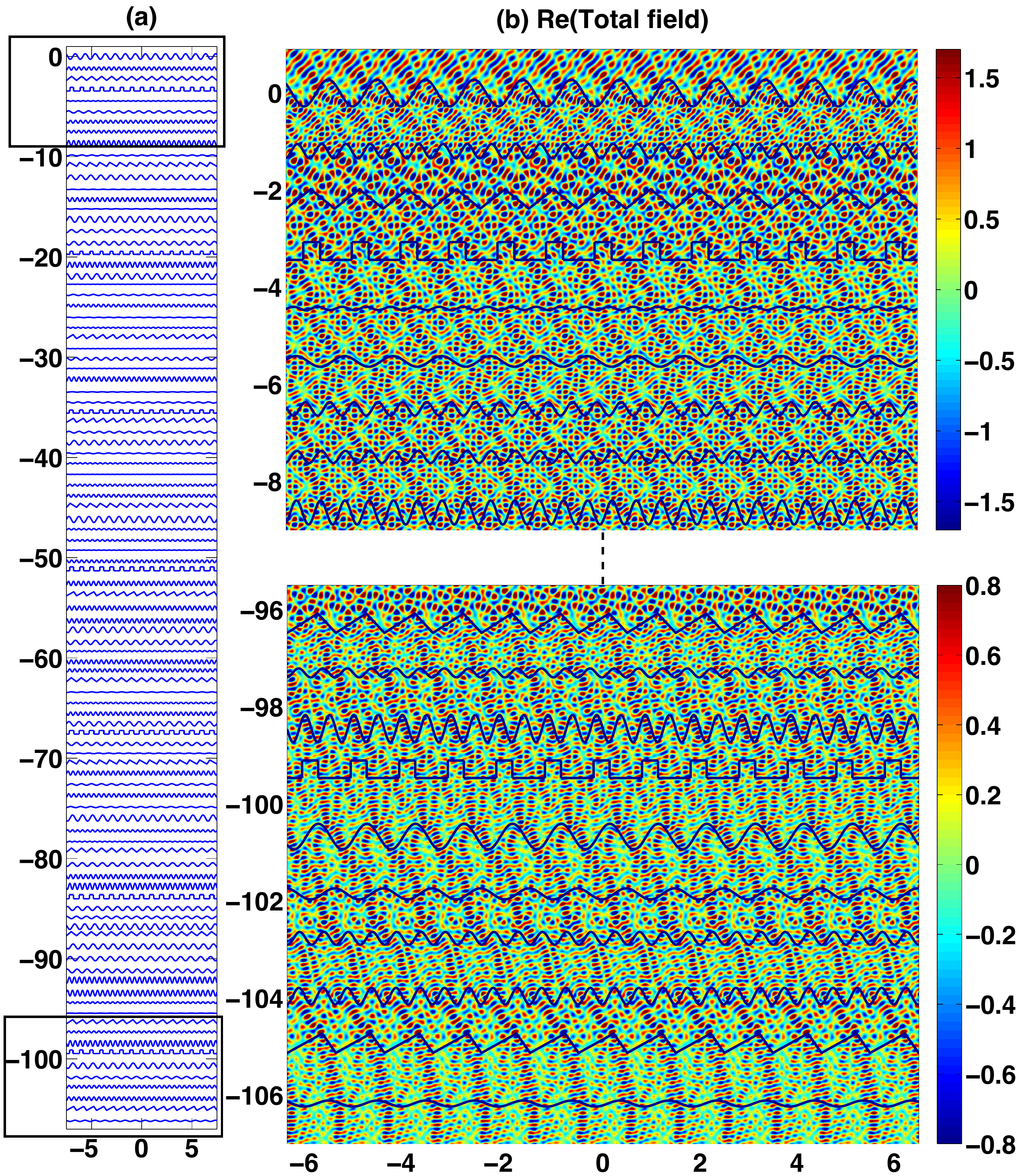} 
   \caption{{\bf (a)} The 100-interface structure tested at a Wood anomaly for the top layer. {\bf (b)} Real part of total field $u+\ui$ in the rectangles drawn in (a), for 
$\omega = 9\pi$, $\ti = -\cos^{-1}(1-2\pi/\omega)$,
$\varepsilon_1 = 1$ and all other $\varepsilon_i$ are randomly chosen between 1 and 2.
$N_i=260$ on sines, $100\times 2$ on triangles, $90\times 5$ on rectangles, $M_w=120$,
$M=60$, $P=120$, and $K=10$.
Flux error is $5\times 10^{-10}$, total solution time 35 sec (not including field evaluation).
(Color online.) }
   \label{field_shape100}
\end{figure}

\begin{figure}[ht] 
   \centering
   \includegraphics[width=5.0in]{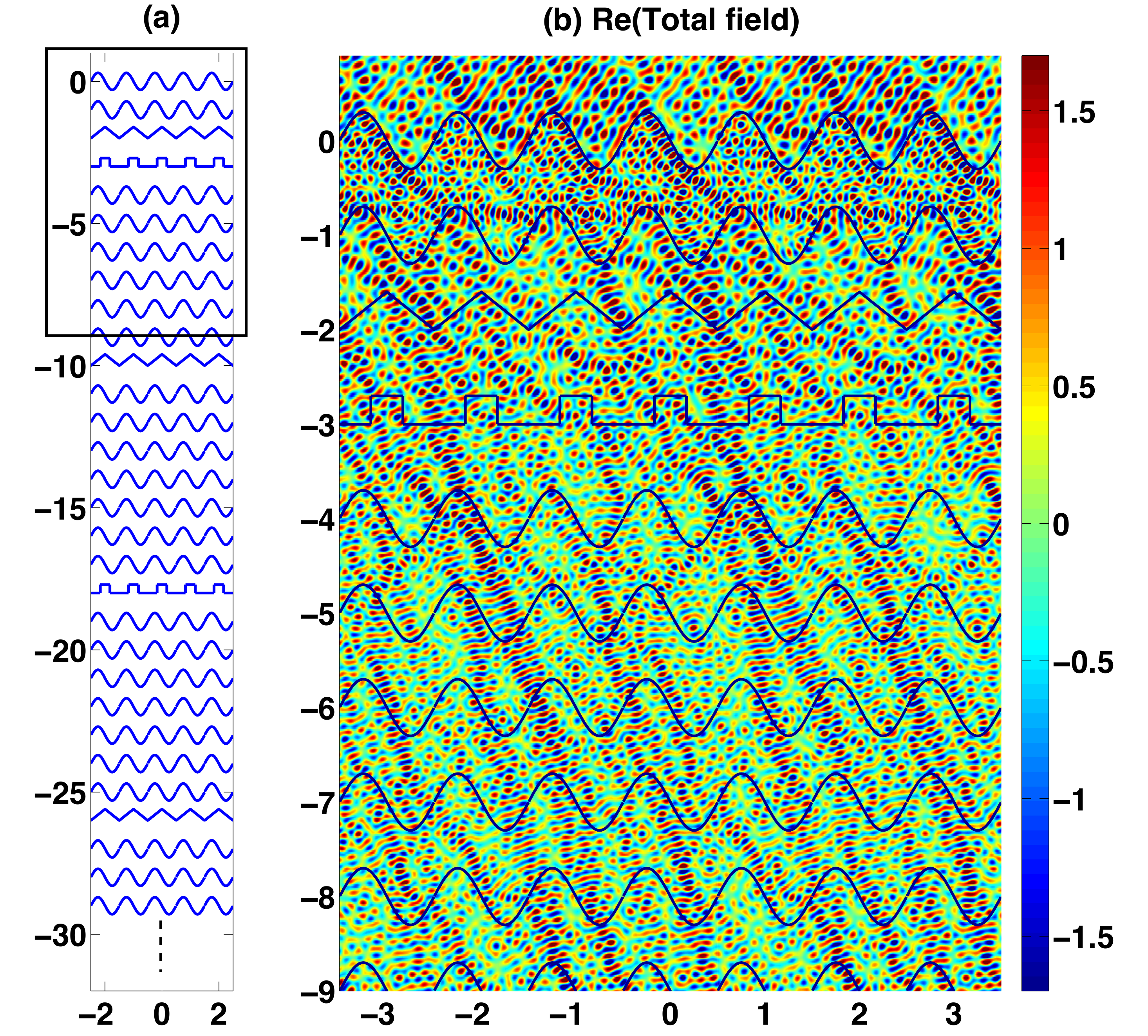} 
   \caption{{\bf (a)} First 30 interfaces of the 1000-interface structure
used in Tables~\ref{table_omega_5} and \ref{table_omega_40}. {\bf (b)} Real part of the total field $u+\ui$ in the rectangle drawn in (a): $\omega = 40$, $\theta^{inc} = -\pi/5$, $\varepsilon_1 = 1$, and all other $\varepsilon_i$ randomly chosen between 1 and 2. Flux error is $4.7\times 10^{-8}$. (Color online.) }
   \label{field_shape1000}
\end{figure}

\begin{figure}[ht] 
   \centering
   \includegraphics[width=4.6in]{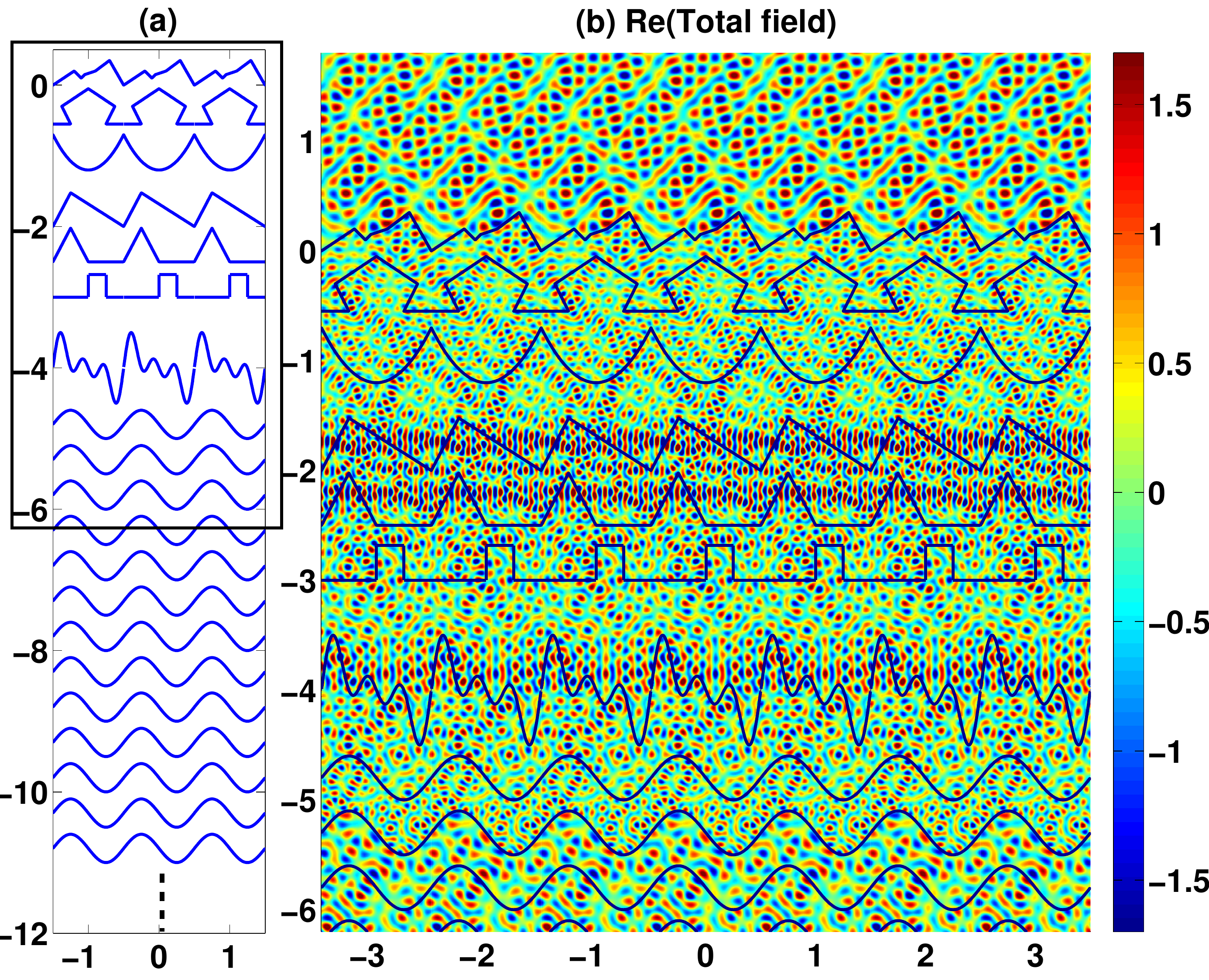} 
   \caption{{\bf (a)} 1000-interface structure consisting of 7 complex-shaped interfaces on top of 993 sine interfaces. {\bf (b)} Real part of total field $u+\ui$ in the rectangular region drawn in (a): $\omega = 40$, $\theta^{inc} = -\pi/4$, $\varepsilon_1 = 1$ and all other $\varepsilon_i$ are chosen randomly between 1 and 3. $N_1 = 160 \times 6$, $N_2 = 160 \times 6$, $N_3=250$,  $N_4=180 \times 2$,  $N_5=160\times 3$, $N_6 = 160\times 5$, and $N_7=300$ on the first 7 interfaces. $N_i = 300$ for the rest of the sine interfaces. $M_w$ = 130, $M$ = 80, $P = 150$, $K$ = 20, $R = 1.5$. Flux error is $7\times 10^{-9}$, total matrix filling time 192 sec, Schur complement 107 sec, block matrix solve 103 sec, total memory used 28 GB. Field evaluation ($1000 \times 1000$ grid points) took 446 sec.
(Color online.) }
   \label{field_general_shape1000}
\end{figure}

\begin{figure}[ht] 
   \centering
   \includegraphics[width=6.0in]{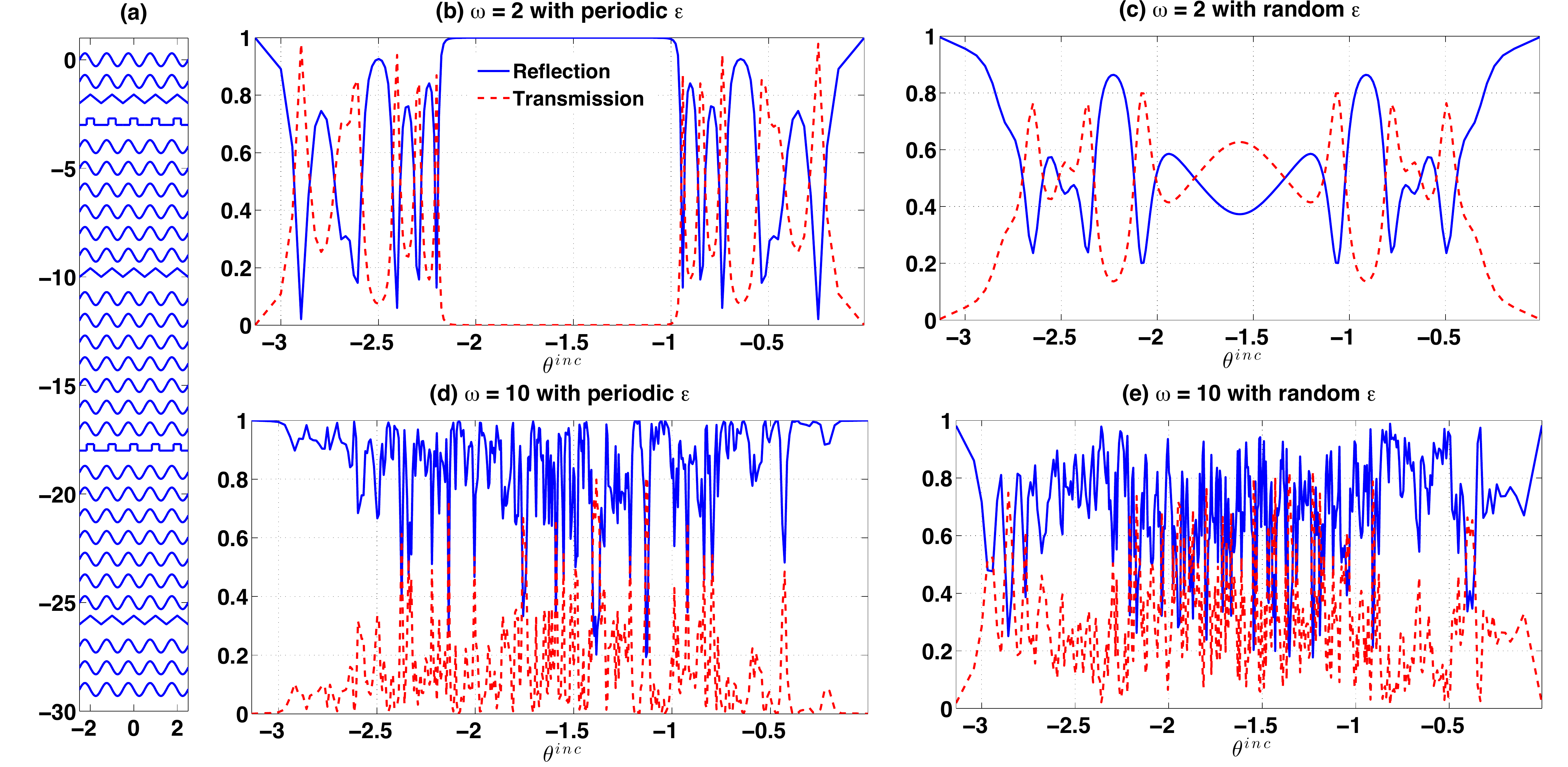} 
   \caption{{\bf (a)} 30-interface structure. Reflection (blue solid line) and transmission (red dashed line) as a function of incident angle from $-\pi$ to $0$ for: {\bf (b)} $\omega = 2$ with periodic $\varepsilon = \{1,4,1,4,\cdots, 4, 1, 4, 1\}$, average flux error $8.3\times 10^{-9}$; {\bf (c)} $\omega=2$ with $\varepsilon_1 = 1$ and all other $\varepsilon_i$ random between 1 and 4, average flux error $1.3\times 10^{-10}$; {\bf (d)}
same structure as (b) but $\omega = 10$,
average flux error $4.7\times 10^{-7}$;
 and {\bf (e)} same structure as (c) but $\omega = 10$,
average flux error $4.1\times 10^{-10}$. (Color online.)}
   \label{reflection_transmission}
\end{figure}

Secondly, in order to study convergence for non-smooth interfaces, every other sine interface is replaced by a rectangular-ridge interface consisting of five line segments
(see the inset in Fig.~\ref{convergence_rectangle}(c)),
with Kress grading parameter $q=6$.
We study the convergence of the two types of interfaces independently.
Fig.~\ref{convergence_rectangle}(a) shows convergence in $N$ on sine interfaces, where $N_i = N$, $i=1,3,\cdots, 29$, while the quadrature points on rectangle is fixed at $N_i = 110 \times 5$, $i = 2,4,\cdots, 30$ 
(i.e.\ $n_{i,l}=110$ for all $l=1,\ldots,5$).
It shows the same convergence as before. Then, $N_i$ is fixed at 70 on all sine interfaces ($i=1,3, \cdots, 29$), and $N_i = N \times 5$ on all rectangle interfaces ($i=2,4,\cdots, 30$), is increased. Both the scattered field and flux error appear to converge as
$\bigO(N^{-6})$ in Fig.~\ref{convergence_rectangle}(b), the order expected from the $q=6$
grading.
Fig.~\ref{convergence_rectangle}(c) shows super-algebraic convergence in $P$.

These tests confirm that flux error is a good indicator of pointwise
error in $u$; from now on we quote flux error.

\subsection{Performance}

Tables \ref{table_omega_5} and \ref{table_omega_40} present the CPU time, memory usage, and flux error for periodic dielectric structures with $I=1$,
3, 10, 30, 100, 300, and 1000 mixed sine, triangle, and rectangle interfaces.
The first is at frequency $\omega=5$ (0.8 wavelengths per period), the second
$\omega = 40$ (6.4 wavelengths per period).
Triangle interfaces are placed every 8th interface starting from the 4th interface and rectangle interfaces are placed every 15th interface starting from the 4th interface. All other interfaces are sine.
Fig.~\ref{field_shape1000}(a) shows the top 30 interfaces of the structure.
All other parameters used for computations are specified in the table captions.
In both tables, the computation time to solve the matrix (Schur complement and block matrix solve) and memory usage increases linearly as number of interfaces increases as expected in \sref{fastsolver}. For 1000 interfaces, it took 151 sec and 634 sec for
$\omega = 5$ and $40$, respectively; note that more quadrature nodes are needed
to accurately discretize the more oscillatory functions for higher $\omega$.
The sizes of the full matrix \eqref{linsys} for 1000 interfaces are $468260\times 287922$ and $832000 \times 751762$, respectively.
At $\omega=40$ the structure is around $8000 \lambda$ tall.

In order to present a performance of the numerical method in some extreme cases, three numerical examples are presented. First, we considered 100 interfaces chosen randomly
as sine, triangle, or rectangle type, with
random heights, phases, layer thicknesses (while preventing collisions),
and layer permittivities; see Fig.~\ref{field_shape100}(a).
We chose a frequency corresponding to a period of 4.5 wavelengths in the top layer,
and an incident angle making the top layer precisely at a Wood anomaly.
Due to space limitations, the real part of total field in only the first 10 and last 10 layers (regions enclosed by rectangles in Fig.\ref{field_shape100}(a)) is plotted in Fig.~\ref{field_shape100}(b).
Matrix filling time was 19 sec, the Schur complements 9 sec, and the tridiagonal solve
8 sec. We achieve 9-digit accuracy in flux;
in general we find that the flux error is no worse at Wood anomalies than at other
angles.

The second and third examples are for 1000 interfaces. The real part of total field $u+u^{inc}$ in only the
top 10 layers is displayed due to figure resolution limitations.
\fref{field_shape1000}(b) is the real part of the total field from the 1000-interface
case in the last column of Table \ref{table_omega_40}.
The third example is presented to highlight the geometric flexibility: we
considered 1000 interfaces consisting of seven complex interface shapes on top of 993 sine interfaces (Fig.~\ref{field_general_shape1000}(a)). Here $\varepsilon_i$, $i=2,\ldots,1000$ are chosen randomly between 1 and 3. All other parameters are given in the figure caption. The total field is computed and displayed in Figure \ref{field_general_shape1000}(b). The CPU time for the solve was about 400 sec with $7\times 10^{-9}$ flux error;
evaluation of the solution at 1 million points took a similar time.

\subsection{Transmission and Reflection Spectrum}

We now compute transmission and reflection spectra as a function of incident angle,
$-\pi < \theta^{inc} < 0$, and benchmark the acceleration technique of \sref{sweep}.
The 30-interface structure shown in Fig.~\ref{reflection_transmission}(a) is used,
the same as in Tables \ref{table_omega_5} and \ref{table_omega_40}.

First, $\omega = 2$ (0.3 wavelengths per period) and periodic permittivities
$\varepsilon = \{1,4,1,4,\cdots, 4, 1\}$ are considered. The spectrum clearly shows Bragg mirror or Fabry-Perot characteristics (there are ranges of incident angles that have total reflection), and symmetry, in Fig.~\ref{reflection_transmission}(b).
Because the wavelength is larger than the geometric features, the interface shape does not play an important role in determining the scattering.
However, when $\varepsilon_i$ are set to random numbers between 1 and 4, the total reflection regime disappears in Fig.~\ref{reflection_transmission}(c).
Since $\omega d/\pi<1$, essentially no benefit comes from multiple angles
sharing $\al$ values, but precomputation of matrix blocks does help.
For 200 incident angles, the computation took 50 sec to produce
with acceleration but 352 sec without, a $7\times$ speed-up.

The frequency is now increased to $\omega = 10$ (1.6 wavelengths per period),
with 641 incident angles. Regardless of the $\varepsilon$ distribution, the transmission-reflection spectrum behaves very abruptly in both periodic (Fig.~\ref{reflection_transmission}(d)) and random (Fig.~\ref{reflection_transmission}(e)) $\varepsilon$ cases.
Also notice that the symmetry of the spectrum is broken because the layered structure has rectangle and triangle shaped interfaces which are now {\em resolved} by the
wavelength.
The computation took 121 sec with acceleration,
and 1619 sec without, a speed-up of $13\times$.
Finally, at $\omega = 20$ (3.2 wavelengths per period),
using 1279 incident angles, the speed-up was about $25\times$
(we don't show these spectra since they do not show any new phenomena).
This is consistent with the
acceleration factor growing linearly with $\omega$. Thus for $\omega=40$ as in
\tref{table_omega_40} a speedup of $50\times$ is expected.

\section{Conclusion} \label{conc}

We presented a new robust and fast integral equation method for
2D scattering from a periodic dielectric grating with an arbitrary
number of layers of general shape.
There are three main features: (1) The computational cost of the new method
scales optimally (linearly) in the number of layers,
allowing 1000-layer structures to be solved rapidly.
(2) The method is stable for {\em all} scattering parameters
including Wood anomalies, since it is based on free-space
rather than quasi-periodic Green's functions.
(3) The periodizing scheme is simple, high-order accurate,
and largely supersedes the Sommerfeld integral
method of the second author and Greengard in \cite{qpsc}.
This solver is expected to be useful in variety of wave
applications in engineering and experimental physics,
including the high accuracy modeling and optimization of optical,
electromagnetic, and acoustic devices, and meta-materials.

There are several natural extensions of this work.
Allowing dielectric inclusions in the layers (as in \cite{junlai}), or
material triple-junctions
(e.g.\ incorporating robust representations of Lee--Greengard as in
\cite{triplejuncQPFDS}),
is simply a matter of bookkeeping,
as long as the number of unknowns per layer remains small (e.g.\ less than $10^4$).
Since we use high-order quadrature schemes, this would allow significantly
more complex unit cell shapes than presented here.
Beyond this, an iterative FMM solution
of the combined system would be appropriate when the number of layers is small
\cite{junlai};
for robustness with many layers a hierarchical fast direct solver
such as in \cite{qpfds,triplejuncQPFDS} could be used in each layer, combined
with our tridiagonal block solve.
Our scheme generalizes to 3D without needing new ideas,
given a surface quadrature for the integral operators.
However, the number of unknowns per layer could then easily exceed $10^4$,
demanding something more elaborate than dense direct linear algebra within each layer.

For a production code, matrix filling and evaluation should be
implemented in C or Fortran,
and a parallel implementation would allow more simultaneous filling of matrix
blocks, as well exploiting a parallel tridiagonal solve.
In terms of analysis, an extension of the rigorous framework of free-space
integral equations to include the presented periodizing scheme is needed.

The MATLAB codes which implement the methods of this paper and
generate some of the figures can be downloaded from
{\tt http://math.dartmouth.edu/$\sim$mhcho/software}
\\
These rely on layer-potential quadrature codes in the
{\tt MPSpack} toolbox by the second author, which can be downloaded from
 {\tt http://code.google.com/p/mpspack}

\section*{Acknowledgments}
We benefited from a helpful discussion with Stephen Shipman.
The work of AHB is supported by NSF grant DMS-1216656.

\bibliography{alex}

\begin{thebibliography}{10}

\bibitem{a+s}
M.~Abramowitz and I.~A. Stegun.
\newblock {\em Handbook of Mathematical Functions with Formulas, Graphs, and
  Mathematical Tables}.
\newblock Dover, New York, 10th edition, 1964.

\bibitem{alpert}
B.~K. Alpert.
\newblock Hybrid {G}auss-trapezoidal quadrature rules.
\newblock {\em SIAM J. Sci. Comput.}, 20:1551--1584, 1999.

\bibitem{arenshabil}
T.~Arens.
\newblock Scattering by biperiodic layered media: The integral equation
  approach.
\newblock Habilitation thesis, Karlsruhe, 2010.

\bibitem{atwater}
H.~A. Atwater and A.~Polman.
\newblock Plasmonics for improved photovoltaic devices.
\newblock {\em Nature Materials}, 9(3):205--213, 2010.

\bibitem{pollution}
I.~M. Babuska and S.~A. Sauter.
\newblock Is the pollution effect of the {FEM} avoidable for the {H}elmholtz
  equation considering high wave numbers?
\newblock {\em SIAM J. Numer. Anal.}, 34(6):2392--2423, 1997.

\bibitem{baodobsonrev}
G.~Bao and D.~C. Dobson.
\newblock Modeling and optimal design of diffractive optical structures.
\newblock {\em Surv. Math. Ind.}, 8:37--62, 1998.

\bibitem{mfs}
A.~H. Barnett and T.~Betcke.
\newblock Stability and convergence of the {M}ethod of {F}undamental
  {S}olutions for {H}elmholtz problems on analytic domains.
\newblock {\em J. Comput. Phys.}, 227(14):7003--7026, 2008.

\bibitem{mpspack}
A.~H. Barnett and T.~Betcke.
\newblock {\tt MPSpack}: A {MATLAB} toolbox to solve {H}elmholtz {PDE}, wave
  scattering, and eigenvalue problems, 2008--2014.
\newblock {\tt http://code.google.com/p/mpspack}.

\bibitem{qplp}
A.~H. Barnett and L.~Greengard.
\newblock A new integral representation for quasi-periodic fields and its
  application to two-dimensional band structure calculations.
\newblock {\em J. Comput. Phys.}, 229:6898--6914, 2010.

\bibitem{qpsc}
A.~H. Barnett and L.~Greengard.
\newblock A new integral representation for quasi-periodic scattering problems
  in two dimensions.
\newblock {\em BIT Numer. Math.}, 51:67--90, 2011.

\bibitem{NIF}
C.~P.~J. {Barty, et al.}
\newblock An overview of {LLNL} high-energy short-pulse technology for advanced
  radiography of laser fusion experiments.
\newblock {\em Nuclear Fusion}, 44(12):S266, 2004.

\bibitem{Bo85}
A.~Bogomolny.
\newblock Fundamental solutions method for elliptic boundary value problems.
\newblock {\em SIAM J. Numer. Anal.}, 22(4):644--669, 1985.

\bibitem{bonnetBDS}
{Bonnet-BenDhia, A.-S.} and F.~Starling.
\newblock Guided waves by electromagnetic gratings and non-uniqueness examples
  for the diffraction problem.
\newblock {\em Math. Meth. Appl. Sci.}, 17:305--338, 1994.

\bibitem{brunohaslam09}
O.~P. Bruno and M.~C. Haslam.
\newblock Efficient high-order evaluation of scattering by periodic surfaces:
  deep gratings, high frequencies, and glancing incidences.
\newblock {\em J. Opt. Soc. Am. A}, 26(3):658--668, 2009.

\bibitem{brunoqp3d}
O.~P. Bruno, S.~Shipman, C.~Turc, and S.~Venakides.
\newblock Efficient evaluation of doubly periodic {G}reen functions in {3D}
  scattering, including {W}ood anomaly frequencies, 2013.
\newblock preprint, {\tt arXiv:1307.1176v1}.

\bibitem{chocai12}
M.~H. Cho and W.~Cai.
\newblock A parallel fast algorithm for computing the {H}elmholtz integral
  operator in {3-D} layered media.
\newblock {\em J. Comput. Phys.}, 231:5910--25, 2012.

\bibitem{coltonkress}
D.~Colton and R.~Kress.
\newblock {\em Inverse acoustic and electromagnetic scattering theory},
  volume~93 of {\em Applied Mathematical Sciences}.
\newblock Springer-Verlag, Berlin, second edition, 1998.

\bibitem{fmm2}
W.~Y. Crutchfield, Z.~Gimbutas, G.~L., J.~Huang, V.~Rokhlin, N.~Yarvin, and
  J.~Zhao.
\newblock {\em Remarks on the implementation of the wideband FMM for the
  Helmholtz equation in two dimensions}, volume 408 of {\em Contemp. Math.},
  pages 99--110.
\newblock Amer. Math. Soc., Providence, RI, 2006.

\bibitem{davisrabin}
P.~J. Davis and P.~Rabinowitz.
\newblock {\em Methods of Numerical Integration}.
\newblock Academic Press, San Diego, 1984.

\bibitem{elschner}
J.~Elschner, R.~Hinder, and G.~Schmidt.
\newblock Finite element solution of conical diffraction problems.
\newblock {\em Adv. Comput. Math.}, 16(2--3):139--156, 2002.

\bibitem{qpfds}
A.~Gillman and A.~Barnett.
\newblock A fast direct solver for quasiperiodic scattering problems.
\newblock {\em J. Comput.\ Phys.}, 248:309--322, 2013.

\bibitem{golubvanloan}
G.~H. Golub and C.~F. Van~Loan.
\newblock {\em Matrix computations}.
\newblock Johns Hopkins Studies in the Mathematical Sciences. Johns Hopkins
  University Press, Baltimore, MD, third edition, 1996.

\bibitem{triplejuncQPFDS}
L.~Greengard, K.~L. Ho, and J.-Y. Lee.
\newblock A fast direct solver for scattering from periodic structures with
  multiple material interfaces in two dimensions.
\newblock {\em J. Comput. Phys.}, 258:738--751, 2014.

\bibitem{gumerov}
N.~A. Gumerov and R.~Duraiswami.
\newblock A method to compute periodic sums.
\newblock {\em J. Comput. Phys.}, 272:307--326, 2014.

\bibitem{hafner90}
C.~Hafner.
\newblock {\em The Generalized Multipole Technique for Computational
  Electromagnetics}.
\newblock Artech House Books, Boston, 1990.

\bibitem{hao}
S.~Hao, A.~H. Barnett, P.~G. Martinsson, and P.~Young.
\newblock High-order accurate {N}ystr{\"o}m discretization of integral
  equations with weakly singular kernels on smooth curves in the plane.
\newblock {\em Adv.\ Comput.\ Math.}, 40(1):245--272, 2014.

\bibitem{qpfdtd}
H.~Holter and H.~Steyskal.
\newblock Some experiences from {FDTD} analysis of infinite and finite
  multi-octave phased arrays.
\newblock {\em IEEE Trans.\ Antennae Prop.}, 50(12):1725--1731, 2002.

\bibitem{horoshenkov}
K.~V. Horoshenkov and S.~N. Chandler-Wilde.
\newblock Efficient calculation of two-dimensional periodic and waveguide
  acoustic {G}reen's functions.
\newblock {\em J. Acoust. Soc. Amer.}, 111:1610--1622, 2002.

\bibitem{jackson}
J.~D. Jackson.
\newblock {\em Classical Electrodynamics}.
\newblock Wiley, 3rd edition, 1998.

\bibitem{jobook}
J.~D. Joannopoulos, S.~G. Johnson, R.~D. Meade, and J.~N. Winn.
\newblock {\em Photonic Crystals: Molding the Flow of Light}.
\newblock Princeton Univ. Press, Princeton, NJ, 2nd edition, 2008.

\bibitem{widebandgrating}
G.~A. Kalinchenko and A.~M. Lerer.
\newblock Wideband all-dielectric diffraction grating on chirped mirror.
\newblock {\em J. Lightwave Technology}, 28(18):2743--49, 2010.

\bibitem{kelzenberg}
M.~D. Kelzenberg, S.~W. Boettcher, J.~A. Petykiewicz, D.~B. Turner-Evans, M.~C.
  Putnam, E.~L. Warren, J.~M. Spurgeon, R.~M. Briggs, N.~S. Lewis, and H.~A.
  Atwater.
\newblock Enhanced absorption and carrier collection in {S}i wire arrays for
  photovoltaic applications.
\newblock {\em Nature Materials}, 9(3):239--244, 2010.

\bibitem{transfer_matrix}
D.~Y.~K. Ko and J.~R. Sambles.
\newblock Scattering matrix method for propagation of radiation in stratified
  media: attenuated total reflection studies of liquid crystals.
\newblock {\em J. Opt. Soc. Am. A}, 5:1863--1866, 1988.

\bibitem{kress91}
R.~Kress.
\newblock Boundary integral equations in time-harmonic acoustic scattering.
\newblock {\em Mathl. Comput. Modelling}, 15:229--243, 1991.

\bibitem{LIE}
R.~Kress.
\newblock {\em Linear Integral Equations}, volume~82 of {\em Appl.\ Math.\
  Sci.}
\newblock Springer, second edition, 1999.

\bibitem{junlai}
J.~Lai, M.~Kobayashi, and A.~H. Barnett.
\newblock A fast solver for the scattering from a layered periodic structure
  with multi-particle inclusions, 2014.
\newblock in preparation.

\bibitem{lechleiter}
A.~Lechleiter and D.-L. Nguyen.
\newblock A trigonometric {G}alerkin method for volume integral equations
  arising in {TM} grating scattering.
\newblock {\em Adv. Comput. Math.}, 40(1):1--25, 2014.

\bibitem{RCWAli4}
L.~Li.
\newblock Formulation and comparison of two recursive matrix algorithms for
  modeling layered diffraction gratings.
\newblock {\em J. Opt. Soc. Am. A}, 13:1024--1035, 1996.

\bibitem{RCWAli}
L.~Li.
\newblock Use of {F}ourier series in the analysis of discontinuous periodic
  structures.
\newblock {\em J. Opt. Soc. Am. A}, 13:1870--1876, 1996.

\bibitem{lintonrev}
C.~M. Linton.
\newblock Lattice sums for the {H}elmholtz equation.
\newblock {\em SIAM Review}, 52(4):603--674, 2010.

\bibitem{linton07}
C.~M. Linton and I.~Thompson.
\newblock Resonant effects in scattering by periodic arrays.
\newblock {\em Wave Motion}, 44:165--175, 2007.

\bibitem{mdirect}
P.~Martinsson and V.~Rokhlin.
\newblock A fast direct solver for boundary integral equations in two
  dimensions.
\newblock {\em J. Comp. Phys.}, 205(1):1--23, 2005.

\bibitem{CWnystrom}
A.~Meier, T.~Arens, S.~N. Chandler-Wilde, and A.~Kirsch.
\newblock A {Nystr\"{o}m} method for a class of integral equations on the real
  line with applications to scattering by diffraction gratings and rough
  surfaces.
\newblock {\em J. Integral Equations Appl.}, 12:281--321, 2000.

\bibitem{scatterometry}
R.~Model, A.~Rathsfeld, H.~Gross, M.~Wurm, and B.~Bodermann.
\newblock A scatterometry inverse problem in optical mask metrology.
\newblock {\em J. Phys.: Conf. Ser.}, 135:012071, 2008.

\bibitem{RCWAmoharam}
M.~G. Moharam and T.~G. Gaylord.
\newblock Rigorous coupled-wave analysis of planar-grating diffraction.
\newblock {\em J. Opt. Soc. Am.}, 71:811--818, 1981.

\bibitem{muller}
C.~M{\"{u}}ller.
\newblock {\em Foundations of the Mathematical Theory of Electromagnetic
  Waves}.
\newblock Springer-Verlag, Berlin, New York, 1969.

\bibitem{nedelec-starling}
J.~C. {N\'{e}d\'{e}lec} and F.~Starling.
\newblock Integral equation methods in a quasi-periodic diffraction problem for
  the time-harmonic {M}axwell's equations.
\newblock {\em SIAM J. Math. Anal.}, 22(6):1679--1701, 1991.

\bibitem{nicholas}
M.~J. Nicholas.
\newblock A higher order numerical method for {3-D} doubly periodic
  electromagnetic scattering problems.
\newblock {\em Commun. Math. Sci.}, 6(3):669--694, 2008.

\bibitem{otani08}
Y.~Otani and N.~Nishimura.
\newblock A periodic {FMM} for {M}axwell's equations in {3D} and its
  applications to problems related to photonic crystals.
\newblock {\em J. Comput. Phys.}, 227:4630--52, 2008.

\bibitem{perryMLD}
M.~D. Perry, R.~D. Boyd, J.~A. Britten, D.~Decker, B.~W. Shore, C.~Shannon, and
  E.~Shults.
\newblock High-efficiency multilayer dielectric diffraction gratings.
\newblock {\em Opt.\ Lett.}, 20(8):940--942, 1995.

\bibitem{rokh83}
V.~Rokhlin.
\newblock Solution of acoustic scattering problems by means of second kind
  integral equations.
\newblock {\em Wave Motion}, 5:257--272, 1983.

\bibitem{absorber}
N.~P. Sergeant, M.~Agrawal, and P.~Peumans.
\newblock High performance solar-selective absorbers using coated
  sub-wavelength gratings.
\newblock {\em Opt. Express}, 18(6):5525--5540, 2010.

\bibitem{shipmanreview}
S.~Shipman.
\newblock {\em Resonant scattering by open periodic waveguides}, volume~1 of
  {\em Progress in Computational Physics (PiCP)}, pages 7--50.
\newblock Bentham Science Publishers, 2010.

\bibitem{stratton}
J.~A. Stratton.
\newblock {\em Electromagnetic Theory}.
\newblock John Wiley \& Sons, Hoboken, NJ, 2007.

\bibitem{taflove}
A.~Taflove.
\newblock {\em Computational Electrodynamics: The Finite-Difference Time-Domain
  Method}.
\newblock Artech House, Norwood, MA, 1995.

\bibitem{refmodels}
G.~L. Wojcik, J.~{Mould Jr.}, E.~Marx, and M.~P. Davidson.
\newblock Numerical reference models for optical metrology simulation.
\newblock In {\em SPIE Microlithography ‘92: IC Metrology, Inspection, and
  Process Control VI, \#1673-06}, 1992.

\bibitem{wood}
R.~W. Wood.
\newblock On a remarkable case of uneven distribution of light in a diffraction
  grating spectrum.
\newblock {\em Philos. Mag.}, 4:396­--408, 1902.

\bibitem{zhaodet}
L.~Zhao and A.~H. Barnett.
\newblock Robust and efficient solution of the drum problem via {Nystr\"om}
  approximation of the {F}redholm determinant, 2014.
\newblock {\tt arxiv:1406.5252}, submitted to {\it J. Comput. Phys.}

\end{thebibliography}

\end{document}